\documentclass[10pt]{article}
\usepackage{lmodern}
\usepackage{amsmath}
\usepackage[T1]{fontenc}
\usepackage[utf8]{inputenc}
\usepackage{authblk}
\usepackage{amsfonts}
\usepackage{graphicx}
\usepackage{rotating}
\usepackage{amssymb}
\usepackage[english]{babel}
\usepackage{xcolor}
\usepackage{amsthm}
\usepackage{graphicx}
\usepackage{mathrsfs}
\usepackage{makecell}
\usepackage{microtype}
\usepackage{mathscinet}
\usepackage{array}
\usepackage{multirow}
\usepackage{enumerate}
\usepackage[cal=boondoxo,bb=ams]{mathalfa}
\usepackage{hyperref}
\hypersetup{hidelinks}

\newtheorem{defn}{Definition}[section]
\newtheorem{theorem}{Theorem}[section]
\newtheorem{prop}{Proposition}[section]
\newtheorem{lemma}{Lemma}[section]

\newtheorem{remark}{Remark}[section]
\newtheorem{exam}{Example}[section]

\newcommand{\ml}{\mathcal}
\newcommand{\mb}{\mathbb}
\DeclareMathOperator{\lin}{lin}
\DeclareMathOperator{\nlin}{nlin}

\def\Xint#1{\mathchoice
	{\XXint\displaystyle\textstyle{#1}}%
	{\XXint\textstyle\scriptstyle{#1}}%
	{\XXint\scriptstyle\scriptscriptstyle{#1}}%
	{\XXint\scriptscriptstyle\scriptscriptstyle{#1}}%
	\!\int}
\def\XXint#1#2#3{{\setbox0=\hbox{$#1{#2#3}{\int}$ }
		\vcenter{\hbox{$#2#3$ }}\kern-.6\wd0}}

\title{On the threshold nature of the Dini continuity for a Glassey derivative-type nonlinearity in a critical semilinear wave equation}
\author[1]{Wenhui Chen\thanks{Wenhui Chen (wenhui.chen.math@gmail.com)}}
\affil[1]{School of Mathematics and Information Science, Guangzhou University,\authorcr 510006 Guangzhou, China}
\author[2]{Alessandro Palmieri\thanks{Alessandro Palmieri (alessandro.palmieri@uniba.it)}}
\affil[2]{Department of Mathematics, University of Bari, 70125 Bari, Italy}
\date{}

\begin{document}
		\maketitle

		\begin{abstract}
			\medskip
		In the present manuscript, we determine the critical condition for the nonlinearity in a semilinear wave equation with a derivative-type nonlinearity. More precisely, we consider a nonlinear term depending on the time derivative of the solution, which is the product of a power nonlinearity with critical Glassey exponent and a modulus of continuity. By employing Zhou's approach along a certain characteristic line, we prove the blow-up in finite time for classical solutions (under a suitable sign condition for the Cauchy data) and we derive upper bound estimates for the lifespan for a not Dini continuous modulus of continuity. Furthermore, in the 3-dimensional and radially symmetric case, by using weighted $L^{\infty}$ estimates, we establish the global existence of small data solutions for a Dini continuous modulus of continuity, and lower bound estimates for the lifespan in the not Dini continuous case. These results provide the regularity threshold (i.e. the Dini condition) for the modulus of continuity in the nonlinearity.
			\\
			
			
			\noindent\textbf{AMS Classification (2020)} 35L05, 35L71, 35B33, 35B44, 35A01
		\end{abstract}
\fontsize{12}{15}
\selectfont

\section{Introduction}\label{Section_Introduction}
\hspace{5mm} Over the past forty years, the following Cauchy problem for  semilinear wave equations with \emph{derivative-type nonlinearity} has been deeply studied:
\begin{align}\label{Wave_ut}
	\begin{cases}
		u_{tt}-\Delta u=|u_t|^{p},&x\in\mb{R}^n,\ t>0,\\
		u(0,x)=\varepsilon u_0(x), \ u_t(0,x)=\varepsilon u_1(x),&x\in\mb{R}^n,
	\end{cases}
\end{align}  where the exponent of the semilinear term satisfies $p>1$, and $\varepsilon>0$ is a parameter describing the size of the Cauchy data. Nowadays, it is known that the critical exponent for the semilinear Cauchy problem \eqref{Wave_ut} is given by the so-called \emph{Glassey exponent}
\begin{align*}
p_{\mathrm{Gla}}(n):=\frac{n+1}{n-1}
\end{align*} for $n\geqslant 2$. By the \emph{critical exponent} here we mean a threshold value in the range for the power of the nonlinear term $|u_t|^p$ such that for $1<p< p_{\mathrm{Gla}}(n)$ (\emph{sub-critical case}) local in time solutions blow up in finite time under suitable sign assumptions for the Cauchy data and regardless of their size, while for $p>p_{\mathrm{Gla}}(n)$  (\emph{super-critical case}) solutions are globally in time defined (in a suitable function space) provided that the Cauchy data are sufficiently small. For this classical model \eqref{Wave_ut}, the \emph{critical case} $p=p_{\mathrm{Gla}}(n)$ belongs to the blow-up range as well. The exponent $p_{\mathrm{Gla}}(n)$ is named after Robert T. Glassey, who conjectured the exact value of the critical exponent (in the sense we just explained) for \eqref{Wave_ut} in \cite{Glassey=1983}, the \emph{Mathematical Review} to the paper \cite{Sideris=1983}. Furthermore, since in the 1-dimensional case local solutions to \eqref{Wave_ut} blow up for any $p>1$ (under sign assumptions for the data), it is customary to set $p_{\mathrm{Gla}}(1):=+\infty$.
Moreover, it is known that for $1<p \leqslant p_{\mathrm{Gla}}(n)$ the \emph{lifespan} $T_\varepsilon$ of a local solution (i.e. the maximal existence time for the solution)  satisfies the following sharp estimates:
\begin{align}\label{Sharp-Lifespan-Glassey}
T_{\varepsilon}\approx \begin{cases}
\varepsilon^{-\left(\frac{1}{p-1}-\frac{n-1}{2}\right)^{-1}}&\mbox{if}\ 1<p<p_{\mathrm{Gla}}(n),\\
\exp\big(C\varepsilon^{-(p-1)}\big)&\mbox{if}\ p=p_{\mathrm{Gla}}(n),
\end{cases}
\end{align} provided that $\varepsilon\in(0,\varepsilon_0]$ for some $\varepsilon_0=\varepsilon_0(n,p,u_0,u_1)>0$.
By writing $T_\varepsilon=+\infty$  if $p>p_{\mathrm{Gla}}(n)$ we can shortly summarize the global existence results.

 Let us outline briefly the main contributions in the proof of the above mentioned \emph{Glassey conjecture} and in the derivation of the sharp lifespan estimates in the sub-critical and critical cases (we are going to emphasize in each case whether the lifespan estimates were obtained from above or from below). The study of \eqref{Wave_ut} was initiated by Fritz John, who proved that for $p=2$ classical solutions to \eqref{Wave_ut} blow up in the 3-dimensional case in \cite{John=1981}, among other results for more general nonlinear terms. Then, in \cite{Sideris=1983} it was proved the global existence of small data and radial solutions in the super-critical case $p>2$ for $n=3$. In \cite{Masuda=1984},  John's blow-up result was extended to space dimensions $n=1,2$. The blow-up results in space dimension $n=2$ were proved for the first time in \cite{Schaeffer=1986,Agemi=1991} for the critical case ($p=3$) and the sub-critical case ($1<p<3$), respectively. For higher space dimensions ($n\geqslant 4$)  blow-up results were proved for radially symmetric solutions in \cite{Rammaha=1987}. The global existence in the super-critical case for $n=2,3$ without radial symmetry was proved independently in \cite{Hidano-Tsutaya=1995,Tzvetkov=1998}. Then, Yi Zhou introduced in \cite{Zhou=2001} a simpler proof of the blow-up results for classical solutions and for any $n\geqslant 1$; moreover, he derived sharp upper bound estimates for the lifespan in the general not radially symmetric case. The global existence in the super-critical case for $n\geqslant 4$ was proved in \cite{Hidano-Wang-Yokoyama=2012} for radially symmetric solutions. Furthermore,  in \cite{Hidano-Wang-Yokoyama=2012} lower bound estimates for the lifespan were derived in the sub-critical ($n\geqslant 2$) and critical ($n\geqslant 3$) cases. The lower bound estimates for the lifespan in the critical case for $n=2$ were derived in \cite{Fang-Wang=2013}. The upper bound estimates for the lifespan and the blow-up results for weak solutions were obtained with a modified test function method in \cite{Ikeda-Sobajima-Wakasa=2019} for all space dimensions. Finally,  quite recently, the lower bound estimates for the lifespan were proved in the 1-dimensional case (as a special unweighted case) in \cite{Kitamura-Morisawa-Takamura=2023}.

In the last decade, the problem considered in the Glassey conjecture has been also investigated for some hyperbolic equations that are different from the classical wave equation  on the Minkowski spacetime. On the one hand, we have results for wave equations on manifolds equipped with non-flat metrics: on asymptotically Euclidean manifolds \cite{Wang=2015,Liu-Wang=2020} and on Friedmann-Lema\^itre-Robertson-Walker spacetimes \cite{Tsutaya-Wakasugi=2021,Hamouda-Hamza-Palmieri=2021,Hamouda-Hamza-Palmieri=2022}. On the other hand, several Glassey-type exponents have been established from the  necessary part in semilinear wave models with time-dependent coefficients: we recall the results for damped wave equations with summable time-dependent coefficients for the linear damping terms \cite{Lai-Takamura=2019,Palmieri-Takamura=2020}, for the wave equation with scale-invariant damping and mass terms \cite{Palmieri-Tu=2021,Hamouda-Hamza=2021}, and for the Tricomi equation \cite{Lucente-Palmieri=2021,Lai-Schiavone=2022,Hamouda-Hamza=2021proc}.

  We can rephrase the critical nature of the Glassey exponent for \eqref{Wave_ut} by saying that $p=p_{\mathrm{Gla}}(n)$ is the threshold value that separates the blow-up range from the small data solutions'  global existence range in the scale of power nonlinearities $\{|u_t|^p:p>1\}$. Following the ideas of \cite{Ebert-Girardi-Reissig=2020}, we wonder what happens if we consider the scale of nonlinear terms $\{|u_t|^{p_{\mathrm{Gla}}(n)}\mu(|u_t|):\mu$ modulus of continuity$\}$, where $\mu:[0,+\infty)\to [0,+\infty)$ is a \emph{modulus of continuity} if $\mu$ is a continuous, concave and increasing function such that $\mu(0)=0$. 
  
   In other words, we consider the following semilinear Cauchy problem for the wave equation with derivative-type nonlinearity:
  \begin{align}\label{Wave_Modulus_ut}
\begin{cases}
u_{tt}-\Delta u=|u_t|^{p_{\mathrm{Gla}}(n)}\mu(|u_t|),&x\in\mb{R}^n,\ t>0,\\
u(0,x)=\varepsilon u_0(x), \ u_t(0,x)=\varepsilon  u_1(x),&x\in\mb{R}^n,
\end{cases}
\end{align}
for $n\geqslant 2$. The additional factor $\mu(|u_t|)$ on the right-hand side of the equation in \eqref{Wave_Modulus_ut} provides an additional regularity to the nonlinear term. In particular, for H\"older's moduli of continuity, i.e., for the class of moduli of continuity $\{\tau^\beta\}_{\beta\in (0,1]}$, the above mentioned results from the literature tell us that any H\"older's modulus of continuity  belongs to the class of moduli of continuity that guarantees the global existence of small data solutions. Our main goal in the present manuscript is to prescribe a regularity threshold for the modulus of continuity $\mu$ in the nonlinearity, i.e. the Dini condition, that separates the global existence of small data solutions from the blow-up of local solutions to \eqref{Wave_Modulus_ut}. In this sense, our aim is to show that the scale $\{|u_t|^p:p>1\}$ is too rough to describe the critical condition for the nonlinear term $|u_t|^{p_{\mathrm{Gla}}(n)}\mu(|u_t|)$. Indeed, we firstly prove that for a modulus of continuity that does not satisfy the Dini condition (cf. Definition \ref{Def Dini}), local in time classical solutions cannot be globally in time prolonged if the moment of order 0 of the second data is positive (see Theorem \ref{Thm-Blow-up-sub}). Clearly, since we need to consider a modulus of continuity $\mu$ which gives a regularity weaker than any H\"older's regularity, we shall consider a \emph{$Log^{-\alpha}$-regularity} with $\alpha>0$, or even a weaker regularity. For the nomenclature 
 of the moduli of continuity we address the reader to the introduction of \cite{Lorenz-Reissig=2019}.  Moreover, in the radially symmetric and 3-dimensional case we prove that the Dini condition is not only necessary but also sufficient for the global existence of small data classical solutions (see Theorem \ref{Thm GESDS rad 3d}).

Next, we shortly review the results for semilinear wave equations having as the nonlinear terms a critical power nonlinearity multiplied by an additional factor involving a modulus of continuity weaker than any H\"older's modulus of continuity. In \cite{Ebert-Girardi-Reissig=2020}, the critical condition for the nonlinearity was investigated for the semilinear damped wave equation. To the knowledge of the authors, \cite{Ebert-Girardi-Reissig=2020} is the first paper where the problem of determining the  critical condition for the nonlinearity is studied for a Cauchy problem associated to a semilinear wave model. More precisely, the nonlinear term $|u|^{p_{\mathrm{Fuj}}(n)}\mu(|u|)$ was considered in \cite{Ebert-Girardi-Reissig=2020} and a sharp integral condition for $\mu$ in a neighborhood of 0 was derived (for $n=1,2$ in the sufficiency part). We stress that for the semilinear damped wave equation the critical exponent for the power nonlinearity is the so-called Fujita exponent $p_{\mathrm{Fuj}}(n)=1+2/n$, and coincides with the critical exponent for the corresponding Cauchy problem for the semilinear heat equation. Remarkably, the threshold integral condition on $\mu$ in \cite{Ebert-Girardi-Reissig=2020} coincides with the one for our model in \eqref{Wave_Modulus_ut}. 

Later, in \cite{Dao-Reissig=2021} these results were generalized to the corresponding weakly coupled system. In \cite{M-Djaouti-Reissig=2023} the problem of determining the  critical condition for the nonlinearity was investigated in the 1-dimensional case for semilinear wave equations either with an effective or a scale-invariant damping term. In some sense, it is not surprising that the threshold condition for $\mu$ is the same one as in \cite{Ebert-Girardi-Reissig=2020} in both cases, since the damping structure in these cases is similar to the one for the classical damping (the critical exponent is $p_{\mathrm{Fuj}}(1)=3$).  Recently, several semilinear evolution models were considered  in \cite{D'Abbicco-Girardi=2023,Girardi=2024} with the nonlinear terms given by the product of a critical power nonlinearity and a modulus of continuity, provided that the critical exponent for the power nonlinearity is obtained by the scaling properties of the partial differential operator. We underline that for the semilinear models studied in  \cite{Girardi=2024} the presence of linear damping terms plays a crucial role in the asymptotic behaviors of the solutions to the corresponding linearized problems.

In \cite{Chen-Reissig=2023} the problem of determining the critical condition for the nonlinearity $|u|^{p_{\mathrm{Str}}(n)}\mu(|u|)$ was considered for the semilinear classical wave equations. Here, $p_{\mathrm{Str}}(n)$ stands for the Strauss exponent (see, for example, the introduction of \cite{Wang-Xu=2012} and references therein for a historical overview on the Strauss conjecture, named after the author of \cite{Strauss=1981}). The exponent $p_{\mathrm{Str}}(n)$ is the positive root of the quadratic equation $(n-1)p^2-(n+1)p-2=0$ (as for the Glassey exponent, in the 1-dimensional case we set $p_{\mathrm{Str}}(1)=+\infty$) and is the critical exponent for the Cauchy problem associated to the semilinear wave equation $u_{tt}-\Delta u=|u|^p$. In particular, by working with an additional $Log^{-\gamma}$-regularity through the modulus of continuity in the nonlinear term $|u|^{p_{\mathrm{Str}}(n)}\mu(|u|)$, in \cite{Chen-Reissig=2023} it was proved that in the sub-critical case $\gamma\in (0,1/p_{\mathrm{Str}}(n))$ energy solutions blow up in finite time (under certain assumptions for the data) and that in the 3-dimensional and radially symmetric  case this threshold is sharp, since for $\gamma>1/p_{\mathrm{Str}}(3)$ the global existence of small data solutions holds in weighted $L^\infty$ spaces. 

In this paper, we are going to determine an integral condition on the modulus of continuity $\mu$, that is the Dini condition, see \eqref{Dini condition} below, whose validity is a criterion to separate the global existence of small data solutions to \eqref{Wave_Modulus_ut} from the blow-up in finite time of  local solutions. We are going to adapt Zhou's approach from \cite{Zhou=2001} to the model in \eqref{Wave_Modulus_ut} with not Dini continuous $\mu$. This method consists in proving the blow-up on a characteristic line for a classical solution and it has been applied to several semilinear hyperbolic models with derivative-type nonlinearity, cf. \cite{Palmieri-Tu=2021, Lucente-Palmieri=2021,Hamouda-Hamza-Palmieri=2021,Palmieri-Takamura=2023}. As byproduct of this technique, we also obtain upper bound estimates for the lifespan for small data. Furthermore, we are going to show  that the threshold condition on $\mu$ is sharp (at least in the radially symmetric and 3-dimensional case), by proving that global in time classical solutions exist provided that the Cauchy data are small and the Dini condition for $\mu$ holds. Indeed, by following the approach of Thomas C. Sideris \cite{Sideris=1983} we are going to prove in the Dini continuous case the global existence in suitable weighted $L^{\infty}$ spaces and to determine lower bound estimates for the lifespan of local solutions in the not Dini continuous case. We stress that the lifespan estimates for radially symmetric solutions in the 3-dimensional case  are sharp, in the sense that the upper and lower bounds for the lifespan $T_\varepsilon$ are the same except for some multiplicative constants independent of $\varepsilon$.

\paragraph{\large Notation} We write $f\lesssim g$ if there exists a positive constant $C$ such that $f\leqslant Cg$. The relation $f\approx g$ holds if and only if $g\lesssim f\lesssim g$. When writing $h(\tau)>0$ [resp. $h(\tau)<0$] as $\tau\to 0^+$, we mean that there exists $\tau_0>0$ such that $h(\tau)>0$ [resp. $h(\tau)<0$]  for any $\tau\in(0,\tau_0]$. Moreover,  $B_R$ denotes the ball around the origin with radius $R$ in the whole space and $\omega_n:=\frac{1}{\Gamma(\frac{n}{2}+1)}\pi^{\frac{n}{2}}$ is the $n$-dimensional measure of the unit ball $B_1$ with Euler's gamma function $\Gamma$. The average of $f$ over $\Omega$ with respect to a measure $\nu$ on $\Omega$ is denoted by
\begin{align*}
\Xint{-}_{\Omega}f\,\mathrm{d}\nu:=\frac{1}{\nu(\Omega)}\int_{\Omega}f\,\mathrm{d}\nu.
\end{align*}
Finally, $[y]_+:=\max\{y,0\}$ is the positive part of $y\in\mathbb{R}$,  and for any $k\in\mathbb{N}^*$ we denote, respectively, by $\log^{[k]}$ and $\exp^{[k]}$ the $k$-times iterated logarithmic and exponential functions, namely,
\begin{align*}
& \log^{[k]} := \begin{cases}
\log & \mbox{if} \ k=1, \\
\log \circ \log^{[k-1]}  & \mbox{if} \ k\geqslant 2, 
\end{cases} 
&& \exp^{[k]} := \begin{cases}
\exp & \mbox{if} \ k=1, \\
\exp \circ \exp^{[k-1]}  & \mbox{if} \ k\geqslant 2.
\end{cases}
\end{align*}

\section{Main results}\label{Section_Main_Result}\setcounter{equation}{0}
\hspace{5mm}Let us begin by recalling the notion of \emph{Dini condition} (or \emph{Dini test}) for a modulus of continuity.
\begin{defn}\label{Def Dini}
Let $\mu:[0,+\infty)\to [0,+\infty)$ be a \emph{modulus of continuity}, that is, $\mu$ is a continuous, increasing, concave function and satisfies $\mu(0)=0$. Then, $\mu$ satisfies the Dini condition if 
\begin{align}\label{Dini condition}
\int_{0}^{\tau_0}\frac{\mu(\tau)}{\tau}\, \mathrm{d}\tau <+\infty \qquad \mbox{for some} \ \tau_0>0.
\end{align} 
\end{defn}

\begin{remark} Notice that if \eqref{Dini condition} is satisfied for some $\tau_0>0$ then it holds for any $\tau_0>0$. Indeed, being $\mu$ continuous, increasing and $\mu(0)=0$, the function $\frac{\mu(\tau)}{\tau}$ is locally summable in any compact set $\mathcal{K}\subset (0,+\infty)$. 
\end{remark}

The condition in \eqref{Dini condition} was implicitly introduced in the study of the convergence of the Fourier series in \cite[Page 101]{Dini=1880} (see also \cite{Zygmund=1988}). We emphasize that, albeit different definitions are possible for the notion of modulus of continuity, the Dini condition appears in different areas of modern mathematical analysis: from elliptic equations \cite{Hartman-Wintner=1955,Burch=1978,Kovats=1997} and integral equations \cite{Chang Lara-Saldana=2023}, to dynamical systems \cite{Fan-Jiang=2001}. We refer to \cite[Section 2, Paragraph A.2]{Andrade-Pellegrino-Pimentel-Teixeira=2022} for a class of non-trivial examples of Dini-moduli of continuity obtained through generalized power series.

As we are going to see in Subsection \ref{Subsection main results GESDS}, the Dini condition for $\mu$ is used in order to prove the global well-posedness for small data solutions to \eqref{Wave_Modulus_ut}. On the other hand, for the blow-up results in Subsection \ref{Subsection main results blowup} we work when \eqref{Dini condition} does not hold, namely, when 
\begin{align}
	\label{Assumption-Thm-Blow-up}
	\int_0^{\tau_0}\frac{\mu(\tau)}{\tau}\,\mathrm{d}\tau=+\infty\qquad \mbox{for some} \ \tau_0>0.
\end{align}

\subsection{Blow-up and upper bound estimates for the lifespan of solutions} \label{Subsection main results blowup}
\hspace{5mm}  We begin by stating the blow-up result for semilinear classical wave equations \eqref{Wave_Modulus_ut}.

\begin{theorem}\label{Thm-Blow-up-sub}
Let $n\geqslant 2$ and let $u_0\in \ml{C}^2_0(\mathbb{R}^n)$, $u_1\in \ml{C}^1_0(\mathbb{R}^n)$ be compactly supported functions with supports contained in $B_R$ for some $R>0$. Furthermore, we assume that $u_1$ satisfies 
\begin{align*}
\int_{\mb{R}^n}u_1(x)\,\mathrm{d}x>0.
\end{align*} Let $\mu:[0,+\infty)\to [0,+\infty)$ be a modulus of continuity fulfilling the following assumptions:
\begin{itemize}
\item[$(A_1)$] the function $g:\tau\in\mb{R}\to |\tau|^{p_{\mathrm{Gla}}(n)}\mu(|\tau|)\in [0,+\infty)$ is convex;
\item[$(A_2)$] $\mu$ does not satisfy the Dini condition, i.e., \eqref{Assumption-Thm-Blow-up} holds for some (and, then, for any) $\tau_0>0$.

\end{itemize}
If $u\in\ml{C}^2([0,T_{\varepsilon})\times\mb{R}^n)$ is a classical  local (in time)  solution to the semilinear Cauchy problem \eqref{Wave_Modulus_ut}, then $u$ blows up in finite time, i.e., $T_\varepsilon <+\infty$.
\end{theorem}


\begin{exam} \label{Example-01}	
Let $\tau_0>0$ be sufficiently small. The assumptions $(A_1)-(A_2)$ on the modulus of continuity (MoC) from Theorem \ref{Thm-Blow-up-sub} hold for the following functions $\mu(\tau)$ that on a small interval $[0,\tau_0]$ satisfy $\mu(0)=0$ and for $\tau\in (0,\tau_0]$:
\begin{description}
	\item[(1MoC)] $\mu(\tau)=(\log\frac{1}{\tau})^{-\gamma}$ with $\gamma\in(0,1]$;
	\item[(2MoC)] $\displaystyle{\mu(\tau)=\prod_{j=1}^{k-1}\left(\log^{[j]}\tfrac{1}{\tau}\right)^{-1}\left(\log^{[k]}\tfrac{1}{\tau}\right)^{-\gamma}  }$ with $\gamma\in(0,1]$ and $k\in \mathbb{N}$, $k\geqslant 2$;
	\item[(3MoC)]  $\displaystyle{\mu(\tau)=\prod_{j=1}^{k}\left(\log^{[j]}\tfrac{1}{\tau}\right)^{-\gamma_j}} $ with $k\in \mathbb{N}$, $k\geqslant 2$ and $\{\gamma_j\}_{1\leqslant j\leqslant k}\subset \mathbb{R}$ fulfilling the next conditions $($clearly, we are excluding the trivial case $\gamma_j=0$ for any $j\in\{1,\ldots,k\})$
	\begin{enumerate}[a$)$]
		\item $\gamma_{j_0}>0$ where $j_0:=\min\{1\leqslant j\leqslant k: \gamma_j\neq 0\}$;
		\item either $\gamma_j=1$ for any $j\in\{1,\ldots,k\}$ or $\gamma_{j_1}<1$, where $j_1:=\min\{1\leqslant  j\leqslant k: \gamma_j\neq 1\}$ $($of course, in this second occurrence $\gamma_j\neq 1$ for some $j\in\{1,\ldots,k\})$.
	\end{enumerate}
\end{description}
\end{exam}

\begin{remark} \label{remark expalnation conditions mu}
 We consider the most general case in Example \ref{Example-01}, i.e., $(3\mathrm{MoC})$, and we explain why $\mu$ is a modulus of continuity that satisfies $(A_1)-(A_2)$. The continuity in $\tau=0$ follows from the fact that $\gamma_{j_0}>0$. For the monotonicity and the concavity of $\mu$ in a right neighborhood of $0$ we need to calculate $\mu'(\tau)$ and $\mu''(\tau)$. Since
 \begin{align*}
\frac{\mathrm{d}}{\mathrm{d}\tau} \left(\log^{[h]} \frac{1}{\tau}\right)=- \frac{1}{\tau} \prod_{j=1}^{h-1} \left(\log^{[j]} \frac{1}{\tau}\right)^{-1}, 	
 \end{align*}
 by Leibniz rule we find
\begin{align}
\mu'(\tau)& = \frac{\mu(\tau)}{\tau} \sum_{h=1}^{k} \gamma_h \prod_{j=1}^{h} \left(\log^{[j]} \frac{1}{\tau}\right)^{-1} ,  \label{mu prime}\\
\mu''(\tau)& = \frac{\mu(\tau)}{\tau^2} \bigg[ -\sum_{h=1}^{k} \gamma_h \prod_{j=1}^{h} \left(\log^{[j]} \frac{1}{\tau}\right)^{-1} + \sum_{h=1}^{k} \gamma_h \sum_{\ell=1}^{h} (\gamma_\ell+1) \prod_{j=1}^{\ell} \left(\log^{[j]} \frac{1}{\tau}\right)^{-2} \prod_{j=\ell+1}^{h} \left(\log^{[j]} \frac{1}{\tau}\right)^{-1} \notag \\
& \qquad \qquad + \sum_{h=1}^{k} \gamma_h \prod_{j=1}^{h} \left(\log^{[j]} \frac{1}{\tau}\right)^{-2}  \sum_{\ell=h+1}^{k} \gamma_\ell \prod_{j=h+1}^{\ell} \left(\log^{[j]} \frac{1}{\tau}\right)^{-1}\bigg].  \label{mu second}
\end{align}
Since $\gamma_{j_0}>0$ it follows that $\mu'(\tau)>0$ and $\mu''(\tau)<0$ as $\tau\to 0^+$. Hence, $\mu$ is strictly increasing and concave in $[0,\tau_0]$, provided that $\tau_0$ is sufficiently small.
We check now that the moduli of continuity in Example \ref{Example-01} fulfill $(A_1)$. By \eqref{mu prime} and \eqref{mu second} we see that 
\begin{align}\label{Cond derivative mu}
\tau^{k} \mu^{(k)}(\tau)=o(\mu(\tau))\qquad  \mbox{as}\  \tau\to0^+\ \mbox{for} \ k=1,2.
\end{align}
Consequently, recalling that $g(\tau)= \tau^{p_{\mathrm{Gla}}(n)} \mu(\tau)$ for $\tau>0$, we obtain
\begin{align*}
g''(\tau) & = \tau^{p_{\mathrm{Gla}}(n)-2}  \left[ p_{\mathrm{Gla}}(n)(p_{\mathrm{Gla}}(n)-1)\mu(\tau)+2p_{\mathrm{Gla}}(n)\tau \mu'(\tau)+\tau^2 \mu''(\tau)\right] \\
& = p_{\mathrm{Gla}}(n)(p_{\mathrm{Gla}}(n)-1)\tau^{p_{\mathrm{Gla}}(n)-2}  \mu(\tau) \left[1+o(1)\right].
\end{align*} Thus, $g''(\tau)>0$ as $\tau\to 0^+$ and $g$ is convex in $[0,\tau_0]$, provided that $\tau_0$ is sufficiently small. Let us show that the moduli of continuity in Example \ref{Example-01} satisfy $(A_2)$.
If $\gamma_j=1$ for any $j\in\{1,\ldots,k\}$, obviously
\begin{align*}
\int_0^{\tau_0} \frac{\mu(\tau)}{\tau}\, \mathrm{d}\tau= \int_0^{\tau_0} \tau^{-1}\prod_{j=1}^k \left(\log^{[j]}\frac{1}{\tau}\right)^{-1} \mathrm{d}\tau= \lim_{\epsilon\to 0^+} \left(\log^{[k+1]}\frac{1}{\epsilon}- \log^{[k+1]}\frac{1}{\tau_0}\right)=+\infty.
\end{align*}
 On the other hand, if  $\gamma_{j_1}<1$, then we can fix $\delta_0>0$ such that $\gamma_{j_1}+\delta_0<1$. Since 
\begin{align*}
\lim_{\tau \to 0^+} \left(\log^{[j_1]}\frac{1}{\tau}\right)^{\delta_0} \prod_{j=j_1+1}^{k} \left(\log^{[j]}\frac{1}{\tau}\right)^{-\gamma_j} =+\infty,
\end{align*} being $\tau_0$ sufficiently small, it results
\begin{align*}
\int_0^{\tau_0} \frac{\mu(\tau)}{\tau}\, \mathrm{d}\tau& \gtrsim  \int_0^{\tau_0} \tau^{-1}\prod_{j=1}^{j_1-1} \left(\log^{[j]}\frac{1}{\tau}\right)^{-1} \left(\log^{[j_1]}\frac{1}{\tau}\right)^{-\gamma_{j_1}-\delta_0}  \mathrm{d}\tau \\ &= \frac{1}{1-\gamma_{j_1}-\delta_0}\lim_{\epsilon\to 0^+} \left(\left(\log^{[j_1]}\frac{1}{\epsilon}\right)^{1-\gamma_{j_1}-\delta_0}- \left(\log^{[j_1]}\frac{1}{\tau_0}\right)^{1-\gamma_{j_1}-\delta_0}\right)=+\infty,
\end{align*}
so $\frac{\mu(\tau)}{\tau}\notin L^1([0,\tau_0])$. Finally, we point out that for a modulus of continuity $\mu$ it is crucial the behavior on a right neighborhood of $0$ $($for this reason in the previous considerations we focused on the interval $[0,\tau_0])$ and that $\mu(\tau)$ can be continued to $\tau\in[0,+\infty)$ in such a way that $\mu$ is a continuous, concave and increasing function and the function $g$ $($that provides the nonlinear term$)$ is a convex function. Moreover, since we prolong $g$ to an increasing and convex function in $[0,+\infty)$, by even reflection we obtain that $g$ is convex in $\mathbb{R}$. As we will see in the proof of Theorem \ref{Thm-Blow-up-sub}, the behavior of $\mu(\tau)$ $[$resp. $g(\tau)]$ for large $\tau$ is completely irrelevant while proving the blow-up of a local solution to \eqref{Wave_Modulus_ut}.
\end{remark}	
	
\begin{exam}
 As a counterexample for the condition \eqref{Assumption-Thm-Blow-up} in Theorem \ref{Thm-Blow-up-sub} we consider the family of H\"older's moduli of continuity, that is,  $\{\mu_\beta\}_{\beta \in (0,1]}$ such that  $\mu_{\beta}(\tau)=\tau^{\beta}$ for any $\tau\in [0,+\infty)$ and $\beta\in(0,1]$. Therefore, this condition is consistent with the global (in time) existence results in \cite{Sideris=1983,Hidano-Tsutaya=1995,Tzvetkov=1998,Hidano-Wang-Yokoyama=2012} for the classical semilinear wave equations  with derivative-type nonlinearity $|u_t|^{p_{\mathrm{Gla}}(n)+\beta}$, which correspond to the super-critical case $p>p_{\mathrm{Gla}}(n)$ in \eqref{Wave_ut}.
\end{exam}

\begin{remark}
	Concerning the semilinear wave equation \eqref{Wave_ut} in the critical case  $p=p_{\mathrm{Gla}}(n)$, the previous works \cite{John=1981,Masuda=1984,Schaeffer=1986,Rammaha=1987,Agemi=1991,Zhou=2001} showed that the solution blows up in finite time. By taking $\mu(|u_t|)$ as an additional factor  that multiplies the critical power nonlinearity $|u_t|^{p_{\mathrm{Gla}}(n)}$, for a modulus of continuity $\mu$ fulfilling \eqref{Assumption-Thm-Blow-up}, Theorem \ref{Thm-Blow-up-sub} tells us that a classical solution still blows up in finite time. Notice that we have to consider moduli of continuity weaker than any H\"older's modulus of continuity, therefore, we look for an infinitesimal function $\mu(\tau)$ as $\tau\to 0^+$ but that tends to $0$ slower than any power $\{\tau^\beta\}_{\beta\in(0,1]}$.
\end{remark}

 As byproduct of our approach in proving the blow-up result, we derive an upper bound estimate for the lifespan $T_{\varepsilon}$ of a local solution considered in Theorem \ref{Thm-Blow-up-sub}.

\begin{prop}\label{Prop-Lifespan}
 Let $n\geqslant 2$. We require that the initial data $u_0,u_1$ and the modulus of continuity $\mu$ satisfy the same assumptions as in Theorem \ref{Thm-Blow-up-sub}. Let us introduce the decreasing function 
\begin{align} \label{def H function lifespan}
\ml{H}:  \tau\in (0,\tau_0]\to   \int_{\tau}^{\tau_0}\frac{\mu(\sigma)}{\sigma}\,\mathrm{d}\sigma\in [0,+\infty).
\end{align}
 Then, there exists $\varepsilon_0=\varepsilon_0(n,R,u_1)>0$ such that for any $\varepsilon\in (0,\varepsilon_0]$ the following upper bound estimate for the lifespan of the local (in time) classical solution $u$ holds:
\begin{align}\label{Upper-Lifespan}
T_{\varepsilon}\leqslant C\varepsilon^{\frac{2}{n-1}}\left[\ml{H}^{-1}\left(c_1\varepsilon^{-\frac{2}{n-1}}+\ml{H}(c_2\varepsilon)\right)\right]^{-\frac{2}{n-1}},
\end{align}
where the positive constants $c_1,c_2,C$ are independent of $\varepsilon$.
\end{prop}

\subsection{Global existence and lower bound estimates for the lifespan of solutions} \label{Subsection main results GESDS}
\hspace{5mm}Next, we study the optimality of Dini condition for our model \eqref{Wave_Modulus_ut} in the 3-dimensional and radially symmetric case: in this case, we will prove that for a modulus of continuity satisfying the Dini condition (and suitable regularity and growth conditions) there exists a uniquely determined solution provided that the Cauchy data are sufficiently small. Furthermore, we will show for moduli of continuity not Dini continuous that the upper bound estimate for the lifespan in \eqref{Upper-Lifespan} is sharp, by establishing the corresponding lower bound estimate.

Before stating the global (when $\mu$ is Dini continuous) and the local (when $\mu$ is not Dini continuous) well-posedness for $n=3$, by denoting $r=|x|$ hereafter, let us rewrite the radially symmetric version of \eqref{Wave_Modulus_ut} as
\begin{align}\label{Wave_Modulus_ut_rad 3d}
\begin{cases}
\displaystyle{u_{tt}-u_{rr}-\frac{2}{r}u_r=|u_t|^{2}\mu(|u_t|),} & r\geqslant 0,\ t>0,\\[0.5em]
u(0,r)=\varepsilon f_0(r), \ u_t(0,r)=\varepsilon  f_1(r) ,&r\geqslant 0,
\end{cases}
\end{align} where we take the radial data $u_0(x)=f_0(|x|)$, $u_1(x)=f_1(|x|)$. 

In the next theorem, we study the global (in time) well-posedness of \eqref{Wave_Modulus_ut_rad 3d}. 

\begin{theorem}\label{Thm GESDS rad 3d}
Let $n=3$ and let $f_0\in \mathcal{C}^2([0,+\infty))$, $f_1\in \mathcal{C}^1([0,+\infty))$ be compactly supported functions with supports contained in $[0,R]$ for some $R>0$ with $f'_0(0)=0$. \\ Let $\mu:[0,+\infty)\to [0,+\infty)$ be a modulus of continuity fulfilling the following assumptions:
\begin{itemize}
\item[$(A_3)$] $\mu \in \mathcal{C}^2((0,+\infty))$ and there exists $\tau_1>0$ such that
\begin{align}\label{mu derivatives growth cond}
|\mu^{(k)}(\tau)|\lesssim \tau^{-k}\mu(\tau)
\end{align} for any $\tau\in (0,\tau_1)$ and $k=1,2$;
\item[$(A_4)$] $\mu$ satisfies the Dini condition, i.e., \eqref{Dini condition} holds for some (and, then, for any) $\tau_0>0$.
\end{itemize}
Then, there exists $\varepsilon_0=\varepsilon_0(f_0,f_1,R,\mu)>0$ such that for any $\varepsilon\in (0,\varepsilon_0] $ the semilinear Cauchy problem \eqref{Wave_Modulus_ut_rad 3d} has a uniquely determined global (in time) classical solution $u=u(t,r)$ such that $u\in\mathcal{C}^1([0,+\infty)^2)$ and $ru(t,r)\in \mathcal{C}^2([0,+\infty)^2)$ subject to the support condition
\begin{align*}
\mathrm{supp}\, u \subset \{(t,r)\in [0,+\infty)^2: r\leqslant t+R \}
\end{align*}  and satisfying for any $t,r\geqslant 0$ the following estimates: 
\begin{align}
|\partial^\alpha u(t,r)|& \lesssim \varepsilon \|(f_0,f'_0,f''_0,f_1,f'_1)\|_{L^\infty([0,+\infty))} (t-r+2R)^{-1} & \forall \alpha \in\mathbb{N}^2: |\alpha|=0,1,\label{Est-L-infty-01}\\
|\partial^\alpha (ru(t,r))|& \lesssim \varepsilon \|(f_0,f'_0,f''_0,f_1,f'_1)\|_{L^\infty([0,+\infty))} (t-r+2R)^{-1} & \forall \alpha \in\mathbb{N}^2: |\alpha|=1,2,\label{Est-L-infty-02}
\end{align} 
where $\partial=(\partial_t,\partial_r)$.
\end{theorem}

\begin{remark} Combining the results from Theorems \ref{Thm-Blow-up-sub} and \ref{Thm GESDS rad 3d}, we see that, at least in the radially symmetric and 3-d case, the Dini condition in \eqref{Dini condition} for the modulus of continuity $\mu$ in \eqref{Wave_Modulus_ut}  is a threshold condition in the sense explained in the introduction.
\end{remark}

\begin{exam}
Let us consider $\tau_0>0$, $k\in\mathbb{N}$, $k\geqslant 1$ and $\{\gamma_j\}_{1\leqslant j\leqslant k}$ such that $\gamma_j\neq 0$ for some $j$. An example of modulus of continuity that satisfies the Dini condition in \eqref{Dini condition}  is given by $\mu(\tau)$ such that for $\tau\in(0,\tau_0]$ 
\begin{align*}
	\mu(\tau)=\prod_{j=1}^{k}\left(\log^{[j]}\frac{1}{\tau}\right)^{-\gamma_j}
\end{align*} provided that $\gamma_{j_0}>0$ and $\gamma_{j_1}>1$ (here $j_0$ and $j_1$ are defined as in Example \ref{Example-01}). Clearly $\mu$ is continuous in $\tau=0$, strictly increasing and concave in $[0,\tau_1]$ with $\tau_1$ being sufficiently small, and verifies $(A_3)$. The proof of these facts follows from \eqref{Cond derivative mu}, that can be proved exactly as in Remark \ref{remark expalnation conditions mu}. Let us show now that this $\mu$ fulfills \eqref{Dini condition}. We fix $\delta_0>0$ such that $\gamma_{j_1}-\delta_0>1$. Since
\begin{align*}
\lim_{\tau \to 0^+} \left(\log^{[j_1]}\frac{1}{\tau}\right)^{-\delta_0} \prod_{j=j_1+1}^{k} \left(\log^{[j]}\frac{1}{\tau}\right)^{-\gamma_j} =0^+
\end{align*} for a sufficiently small $\tau_0$, it holds
\begin{align*}
\int_0^{\tau_0} \frac{\mu(\tau)}{\tau}\, \mathrm{d}\tau& \lesssim  \int_0^{\tau_0} \tau^{-1}\prod_{j=1}^{j_1-1} \left(\log^{[j]}\frac{1}{\tau}\right)^{-1}  \left(\log^{[j_1]}\frac{1}{\tau}\right)^{-\gamma_{j_1}+\delta_0} \mathrm{d}\tau 
= \frac{1}{\gamma_{j_1}-1-\delta_0} \left(\log^{[j_1]}\frac{1}{\tau_0}\right)^{1-\gamma_{j_1}+\delta_0},
\end{align*} so $\frac{\mu(\tau)}{\tau}\in L^1([0,\tau_0])$. This example is the counterpart of (3MoC) in Example \ref{Example-01}.
\end{exam}

Recalling the function  $\mathcal{H}$ introduced in \eqref{def H function lifespan}, we next derive a lower bound estimate for the lifespan $T_{\varepsilon}$ of a local (in time) solution if the Dini condition does not hold, i.e. under the assumption $(A_2)$. 

\begin{prop}\label{Prop lower bound lifespan}
Let $n=3$ and let $f_0\in \mathcal{C}^2([0,+\infty))$, $f_1\in \mathcal{C}^1([0,+\infty))$ be compactly supported functions with supports contained in $[0,R]$ for some $R>0$ with $f'_0(0)=0$. \\ Let $\mu:[0,+\infty)\to [0,+\infty)$ be a modulus of continuity fulfilling the following assumptions:
\begin{itemize}
\item[$(A_2)$]  $\mu$ does not satisfy the Dini condition, i.e., \eqref{Assumption-Thm-Blow-up} holds for some (and, then, for any) $\tau_0>0$;
\item[$(A_3)$] $\mu \in \mathcal{C}^2((0,+\infty))$ and there exists $\tau_1>0$ such that
\eqref{mu derivatives growth cond} holds for any $\tau\in (0,\tau_1)$ and $k=1,2$.
\end{itemize}
Then, there exist $\varepsilon_0=\varepsilon_0(f_0,f_1,R,\mu)>0$ and three constants $K,k_1,k_2>0$ depending on $f_0,f_1,R$  such that for any $\varepsilon\in (0,\varepsilon_0] $ and any $T>0$ with
\begin{align}\label{lower bound lifespan 3d rad-T}
T\leqslant K \varepsilon \left[\mathcal{H}^{-1}\left(k_1\varepsilon^{-1}+\mathcal{H}(k_2 \varepsilon)\right)\right]^{-1},
\end{align}
 the semilinear Cauchy problem \eqref{Wave_Modulus_ut_rad 3d}  has a uniquely determined local (in time) classical solution $u=u(t,r)$ such that $u\in\mathcal{C}^1([0,T]\times [0,+\infty))$ and $ru(t,r)\in \mathcal{C}^2([0,T]\times [0,+\infty))$ subject to the support condition $\mathrm{supp}\, u \subset \{(t,r)\in [0,T]\times [0,+\infty): r\leqslant t+R \}$ and satisfying for any $t\in[0,T]$ and $ r\geqslant 0$ the estimates  \eqref{Est-L-infty-01} and \eqref{Est-L-infty-02}.  In other words, the following lower bound estimate for the lifespan of the local (in time) classical solution $u$ holds:
\begin{align}\label{lower bound lifespan 3d rad}
	T_{\varepsilon}\geqslant K \varepsilon \left[\mathcal{H}^{-1}\left(k_1\varepsilon^{-1}+\mathcal{H}(k_2 \varepsilon)\right)\right]^{-1}.
\end{align}
\end{prop}

\begin{remark} Proposition \ref{Prop lower bound lifespan} provides a lower bound estimate for  local (in time) classical solutions to \eqref{Wave_Modulus_ut_rad 3d}. In particular, \eqref{lower bound lifespan 3d rad} tells us that the upper bound estimate in \eqref{Upper-Lifespan} is sharp when $n=3$ in the radially symmetric case. That is to say, for any $\varepsilon\in (0,\varepsilon_0]$ the lifespan $T_{\varepsilon}$ of a solution to \eqref{Wave_Modulus_ut}  satisfies
\begin{align}\label{Lifespan-General-Sharp}
\begin{cases}
\displaystyle{ T_{\varepsilon}\approx \varepsilon^{\frac{2}{n-1}}\left[\ml{H}^{-1}\left(c_1\varepsilon^{-\frac{2}{n-1}}+\ml{H}(c_2\varepsilon)\right)\right]^{-\frac{2}{n-1}}}&\mbox{if}\ \  \displaystyle{\int_0^{\tau_0}\frac{\mu(\tau)}{\tau}\,\mathrm{d}\tau=+\infty,}\\[1em]
T_{\varepsilon}=+\infty&\mbox{if}\ \ \displaystyle{\int_0^{\tau_0}\frac{\mu(\tau)}{\tau}\,\mathrm{d}\tau<+\infty,}
\end{cases}
\end{align}
in the 3-dimensional and radially symmetric case.
\end{remark}
In Section \ref{Section upper bound lifespan}, we rewrite the sharp lifespan estimate \eqref{Lifespan-General-Sharp} in the special cases $(1\mathrm{MoC})$ and $(2\mathrm{MoC})$ from Example \ref{Example-01}.

\section{Blow-up of the classical solutions}\label{Section_Proof}\setcounter{equation}{0}
\hspace{5mm}To begin with the proof, let us introduce the notation $x=(z,w)\in\mb{R}^n$ with $z\in\mb{R}$ and $w\in\mb{R}^{n-1}$ (recall that we are considering the case $n\geqslant 2$). Then, we may define the next auxiliary functions:
\begin{align*}
\ml{U}(t,z)&:=\int_{\mb{R}^{n-1}}u(t,z,w)\,\mathrm{d}w \qquad \mbox{for any} \ (t,z)\in[0,T_{\varepsilon})\times\mb{R},\\
\ml{U}_j(z)&:=\int_{\mb{R}^{n-1}}u_j(z,w)\,\mathrm{d}w \qquad \ \mbox{for any} \ z\in\mb{R} \ \mbox{and}\ j=0,1.
\end{align*}
The function $\ml{U}=\ml{U}(t,z)$ solves the following nonlinear wave equations in 1-d:
\begin{align*}
\begin{cases}
\displaystyle{\ml{U}_{tt}-\ml{U}_{zz}=\int_{\mb{R}^{n-1}}|u_t(t,z,w)|^{p_{\mathrm{Gla}}(n)}\mu(|u_t(t,z,w)|)\,\mathrm{d}w,}&z\in\mb{R},\ t\in(0,T_{\varepsilon}),\\
\ml{U}(0,z)=\varepsilon\ml{U}_0(z),\ \ml{U}_t(0,z)=\varepsilon\ml{U}_1(z),&z\in\mb{R}.
\end{cases}
\end{align*}
By D'Alembert's formula, recalling the notation $g(\tau)=|\tau|^{p_{\mathrm{Gla}}(n)}\mu(|\tau|)$ from $(A_1)$, we obtain the representation
\begin{align*}
\ml{U}(t,z)=\varepsilon\ml{U}^{\lin}(t,z)+\ml{U}^{\nlin}(t,z),
\end{align*}
where 
\begin{align}
\ml{U}^{\lin}(t,z)&:=\frac{1}{2}\big(\ml{U}_0(z+t)+\ml{U}_0(z-t)\big)+\frac{1}{2}\int_{z-t}^{z+t}\ml{U}_1(y)\,\mathrm{d}y,\notag\\
\ml{U}^{\nlin}(t,z)&:=\frac{1}{2}\int_0^t\int_{z-t+s}^{z+t-s}\int_{\mb{R}^{n-1}}g\big(u_t(s,y,w)\big)\,\mathrm{d}w\,\mathrm{d}y\,\mathrm{d}s.\label{U-nlin}
\end{align}

\subsection{Nonlinear integral inequality with the modulus of continuity}
\hspace{5mm} Our next goal is to determine a nonlinear integral inequality for the function $\ml{U}(t,z)$ on the characteristic line $t-z=R$ involving the modulus of continuity. Due to the property of finite speed of propagation, the assumption $\mathrm{supp}\,u_j\subset B_R$ for $j=0,1$ implies that 
 \begin{align*}
 \mathrm{supp}\,u(t,\cdot)\subset B_{R+t} \qquad  \mbox{for any} \ t\in[0,T_{\varepsilon}).
 \end{align*}
Consequently, for $z\in\mathbb{R}$ such that $|z|\geqslant R+t$ we have that the integrand function in the definition of $\mathcal{U}(t,z)$ satisfies $u(t,z,w)=0$ for any $w\in\mathbb{R}^{n-1}$, because $|(z,w)|^2\geqslant (R+t)^2+|w|^2\geqslant (R+t)^2$. Hence,
\begin{align}\label{Supp-U}
\mathrm{supp}\,\ml{U}(t,\cdot)\subset\big(-(R+t),R+t\big)\ \qquad \mbox{for any}\  t\in[0,T_{\varepsilon}).
\end{align}

Since our nonlinear term $g(u_t)$ is local, if we rewrite the support condition with respect to the variable $w\in\mb{R}^{n-1}$, we find
\begin{align*}
\mathrm{supp}\,g\big(u_t(s,y,\cdot)\big)\subset \mathrm{supp}\,u_t(s,y,\cdot)\subset \left\{w\in\mb{R}^{n-1}:\ |w|^2\leqslant (R+s)^2-y^2\right\}
\end{align*}
 for any $s\in[0,t]$ and any $y\in\mb{R}$ such that $|y|\leqslant R+s$. Consequently, for any $y\in\mathbb{R}$ such that $|y|\leqslant R+s$ we have
 \begin{align*}
 \int_{\mb{R}^{n-1}}g\big(u_t(s,y,w)\big)\,\mathrm{d}w&=\int_{B_{\sqrt{(R+s)^2-y^2}}}g\big(u_t(s,y,w)\big)\,\mathrm{d}w\\
 &=\omega_{n-1}\big((R+s)^2-y^2\big)^{\frac{n-1}{2}}\Xint{-}_{B_{\sqrt{(R+s)^2-y^2}}}g\big(u_t(s,y,w)\big)\,\mathrm{d}w.
 \end{align*}
Thanks to the assumption $(A_1)$, which guarantees the convexity of $g$, by Jensen's inequality we get
\begin{align*}
\int_{\mb{R}^{n-1}}g\big(u_t(s,y,w)\big)\,\mathrm{d}w&\geqslant\omega_{n-1}\big((R+s)^2-y^2\big)^{\frac{n-1}{2}}g\left(\Xint{-}_{B_{\sqrt{(R+s)^2-y^2}}}u_t(s,y,w)\,\mathrm{d}w\right)\\
&=\omega_{n-1}\big((R+s)^2-y^2\big)^{\frac{n-1}{2}}g\left(\frac{1}{\omega_{n-1}}\big((R+s)^2-y^2\big)^{-\frac{n-1}{2}}\ml{U}_t(s,y)\right)
\end{align*}
for any $s\in[0,t]$ and any $y\in[z-t+s,z+t-s]$ such that $|y|\leqslant R+s$. From \eqref{U-nlin} we may estimate
\begin{align*}
\ml{U}^{\nlin}(t,z)&\geqslant\frac{\omega_{n-1}}{2}\int_0^t\int_{z-t+s}^{z+t-s}\big[(R+s)^2-y^2\big]_+^{\frac{n-1}{2}}g\left(\frac{1}{\omega_{n-1}}\big[(R+s)^2-y^2\big]_+^{-\frac{n-1}{2}}\ml{U}_t(s,y)\right)\mathrm{d}y\,\mathrm{d}s\\
&=\frac{\omega_{n-1}}{2}\int_{z-t}^{z+t}\int_0^{t-|z-y|}\big[(R+s)^2-y^2\big]_+^{\frac{n-1}{2}}g\left(\frac{1}{\omega_{n-1}}\big[(R+s)^2-y^2\big]_+^{-\frac{n-1}{2}}\ml{U}_t(s,y)\right)\mathrm{d}s\,\mathrm{d}y,
\end{align*}
where we used Fubini's theorem.

Hereafter, we work on the characteristic line $t-z=R$ for $z\geqslant R$ and shrink the domain of integration from the triangular region
\begin{align*}
\left\{(s,y)\in\mb{R}^2:\ 0\leqslant s\leqslant t-|z-y|,\ z-t\leqslant y\leqslant z+t\right\}
\end{align*}
to the parallelogrammatic sub-region
\begin{align*}
\left\{(s,y)\in\mb{R}^2:\ y-R\leqslant s\leqslant y+R,\ R\leqslant y\leqslant z\right\}.
\end{align*}
For this reason, denoting by
\begin{align*}
M(y;R):=\int_{y-R}^{y+R}\big((R+s)^2-y^2\big)^{\frac{n-1}{2}}\mathrm{d}s
\end{align*}
 the measure of the interval $[y-R,y+R]$ with respect to $\big((R+s)^2-y^2\big)^{\frac{n-1}{2}}\mathrm{d}s$ (when $|y|\leqslant R+s$),  it results
\begin{align*}
\ml{U}^{\nlin}(z+R,z)&\geqslant\frac{\omega_{n-1}}{2}\int_R^z\int_{y-R}^{y+R}g\left(\frac{1}{\omega_{n-1}}\big((R+s)^2-y^2\big)^{-\frac{n-1}{2}}\ml{U}_t(s,y)\right)\big((R+s)^2-y^2\big)^{\frac{n-1}{2}}\mathrm{d}s\,\mathrm{d}y\\
&=\frac{\omega_{n-1}}{2}\int_R^zM(y;R)\Xint{-}_{y-R}^{y+R}g\left(\frac{1}{\omega_{n-1}}\big((R+s)^2-y^2\big)^{-\frac{n-1}{2}}\ml{U}_t(s,y)\right)
\big((R+s)^2-y^2\big)^{\frac{n-1}{2}}\mathrm{d}s\,\mathrm{d}y.
\end{align*}
Let us derive upper and lower bound estimates for $M(y;R)$ for $y\geqslant R$. On the one hand,
\begin{align*}
M(y;R)&=\int_{y-R}^{y+R}\big((R+s-y)(R+s+y)\big)^{\frac{n-1}{2}}\mathrm{d}s\\
&\geqslant(2y)^{\frac{n-1}{2}}\int_{y-R}^{y+R}(R+s-y)^{\frac{n-1}{2}}\mathrm{d}s\\
&=\frac{2^{n+1}}{n+1}R^{\frac{n+1}{2}}y^{\frac{n-1}{2}},
\end{align*}
on the other hand,
\begin{align*}
M(y;R)&\leqslant (2R)^{\frac{n-1}{2}}\int_{y-R}^{y+R}(R+s+y)^{\frac{n-1}{2}}\mathrm{d}s\\
&=\frac{2^{n+1}}{n+1}R^{\frac{n-1}{2}}\big((R+y)^{\frac{n+1}{2}}-y^{\frac{n+1}{2}}\big)\\
& = 2^nR^{\frac{n+1}{2}}(y+\theta_0R)^{\frac{n-1}{2}}\\
&\leqslant 2^{\frac{3n-1}{2}}R^{\frac{n+1}{2}}y^{\frac{n-1}{2}},
\end{align*}
for some $\theta_0\in(0,1)$, where  in the second last step we used the mean value theorem. Summarizing the previous two inequalities, we proved that
\begin{align*}
	M(y;R)\approx R^{\frac{n+1}{2}}y^{\frac{n-1}{2}}\ \qquad \mbox{as}\ y\geqslant R.
\end{align*}

By applying Jensen's inequality in the $s$-integral with respect to the measure $\big((R+s)^2-y^2\big)^{\frac{n-1}{2}}\mathrm{d}s$,  for $z\geqslant R$ we obtain
\begin{align*}
\ml{U}^{\nlin}(z+R,z)&\geqslant\frac{2^n \omega_{n-1}}{n+1}R^{\frac{n+1}{2}}\int_R^zy^{\frac{n-1}{2}}\Xint{-}_{y-R}^{y+R}g\left(\frac{1}{\omega_{n-1}}\big((R+s)^2-y^2\big)^{-\frac{n-1}{2}}\ml{U}_t(s,y)\right)
\big((R+s)^2-y^2\big)^{\frac{n-1}{2}}\mathrm{d}s\,\mathrm{d}y\\
&\geqslant\frac{2^n \omega_{n-1}}{n+1}R^{\frac{n+1}{2}}\int_R^zy^{\frac{n-1}{2}}g\left(\frac{1}{\omega_{n-1}} \ \Xint{-}_{y-R}^{y+R} \ml{U}_t(s,y)\,\mathrm{d}s \right)\mathrm{d}y\\
&=\frac{2^n \omega_{n-1}}{n+1}R^{\frac{n+1}{2}}\int_R^zy^{\frac{n-1}{2}}g\left(\frac{1}{\omega_{n-1}M(y;R)}\int_{y-R}^{y+R}\ml{U}_t(s,y)\,\mathrm{d}s\right)\mathrm{d}y.
\end{align*}
Since $g$ is an increasing function on $[0,+\infty)$, being the product of two nonnegative and increasing functions, taking 
\begin{align*}
\widetilde{K}:=\omega_{n-1}^{-1}2^{-\frac{3n-1}{2}}R^{-\frac{n+1}{2}},
\end{align*}
 we conclude that
 \begin{align}
 \ml{U}^{\nlin}(z+R,z)&\geqslant \frac{2^n \omega_{n-1}}{n+1}R^{\frac{n+1}{2}} \int_R^zy^{\frac{n-1}{2}}g\left(\widetilde{K}y^{-\frac{n-1}{2}}\int_{y-R}^{y+R}\ml{U}_t(s,y)\,\mathrm{d}s\right)\mathrm{d}y \notag \\
 &= \underbrace{\frac{2^n \omega_{n-1}}{n+1}R^{\frac{n+1}{2}}}_{=:C_1}\int_R^zy^{\frac{n-1}{2}}g\left(\widetilde{K}y^{-\frac{n-1}{2}}\ml{U}(y+R,y)\right)\mathrm{d}y \label{Nonlinear-Inte-U}
 \end{align}
for $z\geqslant R$, where we used the fundamental theorem of calculus associated with $\ml{U}(y-R,y)=0$ that follows from the support condition \eqref{Supp-U}.
Our assumption on the initial data shows that
\begin{align}
\ml{U}^{\lin}(z+R,z)&=\frac{1}{2}\big(\ml{U}_0(2z+R)+\ml{U}_0(-R)\big)+\frac{1}{2}\int_{-R}^{2z+R}\ml{U}_1(y)\,\mathrm{d}y=\frac{1}{2}\int_{\mb{R}}\ml{U}_1(y)\,\mathrm{d}y\notag\\
&=\frac{1}{2}\int_{\mb{R}^n}u_1(x)\,\mathrm{d}x=:C_0>0,\label{Ulin-Lower}
\end{align}
since $z\geqslant R$ implies that $2z+R\geqslant R$ and $\mathrm{supp}\, \ml{U}_0,\mathrm{supp}\, \ml{U}_1 \subset(-R,R)$ which is a consequence of the assumption $\mathrm{supp}\, u_0 ,\mathrm{supp}\, u_1 \subset B_R$ for the Cauchy data.

 Finally, combining \eqref{Nonlinear-Inte-U} and \eqref{Ulin-Lower}  we obtain the following nonlinear integral inequality:
\begin{align}\label{U-Lower}
\ml{U}(z+R,z)\geqslant C_0\varepsilon +C_1 \int_R^zy^{\frac{n-1}{2}}g\left(\widetilde{K}y^{-\frac{n-1}{2}}\ml{U}(y+R,y)\right)\mathrm{d}y.
\end{align}

\subsection{Proof of Theorem \ref{Thm-Blow-up-sub}}
\hspace{5mm}Let us introduce the $z$-dependent functional
\begin{align*}
\ml{G}(z):=C_0\varepsilon+C_1\int_R^zy^{\frac{n-1}{2}}g\left(\widetilde{K}y^{-\frac{n-1}{2}}\ml{U}(y+R,y)\right)\mathrm{d}y
\end{align*} for $z\geqslant R$.
Clearly, $\ml{G}(R)=C_0\varepsilon$ and
\begin{align}\label{G' >0}
	\ml{G}'(z)=C_1z^{\frac{n-1}{2}}g\left(\widetilde{K}z^{-\frac{n-1}{2}}\ml{U}(z+R,z)\right)>0.
\end{align}
Hence, $\ml{G}=\ml{G}(z)$ is an increasing function, therefore,
\begin{align}\label{G lb epsilon}
\ml{G}(z)\geqslant \ml{G}(R)=C_0\varepsilon.
\end{align}
We remark that from \eqref{U-Lower} we have $\ml{U}(z+R,z)\geqslant \ml{G}(z)$. Consequently, using again the monotonicity of $g$ on the nonnegative semi-axis, we find
\begin{align}
\ml{G}'(z)&\geqslant C_1z^{\frac{n-1}{2}}g\left(\widetilde{K}z^{-\frac{n-1}{2}}\ml{G}(z)\right)\notag\\
&\geqslant C_2 z^{-\frac{n-1}{2}(p_{\mathrm{Gla}}(n)-1)}[\ml{G}(z)]^{p_{\mathrm{Gla}}(n)}\mu \left(\widetilde{K}z^{-\frac{n-1}{2}}\ml{G}(z)\right)\label{Est-G'}
\end{align}
for $z\geqslant R$ with $C_2:=C_1\widetilde{K}^{p_{\mathrm{Gla}}(n)}$. 

Recalling the value of $p_{\mathrm{Gla}}(n)$, the inequalities \eqref{G lb epsilon} and \eqref{Est-G'} yield
\begin{align*}
\frac{\ml{G}'(z)}{[\ml{G}(z)]^{p_{\mathrm{Gla}}(n)}}\geqslant C_2z^{-1}\mu \left(C_0\widetilde{K}\varepsilon z^{-\frac{n-1}{2}}\right).
\end{align*}
Integrating the last inequality over $[R,z]$, one arrives at
\begin{align*}
\frac{1}{p_{\mathrm{Gla}}(n)-1}\left((C_0\varepsilon)^{1-p_{\mathrm{Gla}}(n)}-[\ml{G}(z)]^{1-p_{\mathrm{Gla}}(n)}\right)&=\int_{R}^z\frac{\ml{G}'(y)}{[\ml{G}(y)]^{p_{\mathrm{Gla}}(n)}}\,\mathrm{d}y\\
&\geqslant C_2\int_R^zy^{-1}\mu \left(C_0\widetilde{K}\varepsilon y^{-\frac{n-1}{2}}\right)\mathrm{d}y.
\end{align*}
By using the change of variable $\tau=C_0\widetilde{K}\varepsilon y^{-\frac{n-1}{2}}$, we derive
\begin{align*}
\int_R^zy^{-1}\mu \left(C_0\widetilde{K}\varepsilon y^{-\frac{n-1}{2}}\right)\mathrm{d}y
=\frac{2}{n-1}\int_{C_0\widetilde{K}z^{-\frac{n-1}{2}}\varepsilon}^{C_0\widetilde{K}R^{-\frac{n-1}{2}}\varepsilon}\frac{\mu(\tau)}{\tau}\,\mathrm{d}\tau.
\end{align*}
That is to say,
\begin{align}\label{Final-Ineq}
(C_0\varepsilon)^{-\frac{2}{n-1}}-[\ml{G}(z)]^{-\frac{2}{n-1}}\geqslant\frac{4C_2}{(n-1)^2}\int_{C_0\widetilde{K}z^{-\frac{n-1}{2}}\varepsilon}^{C_0\widetilde{K}R^{-\frac{n-1}{2}}\varepsilon}\frac{\mu(\tau)}{\tau}\,\mathrm{d}\tau.
\end{align}

Let us assume that $u$ is global (in time) defined, then $\ml{U}(z+R,z)$ and, in turn, $\ml{G}(z)$ are defined for any $z\geqslant R$. Nevertheless, this produces a contradiction as $z\to+\infty$, due to the facts that the left-hand side of \eqref{Final-Ineq} is bounded (recall that thanks to \eqref{G' >0} and \eqref{G lb epsilon}, $\ml{G}(z)$ is an increasing function that is bounded from below by a positive constant), while the right-hand side, by the monotone convergence theorem, tends to
\begin{align*}
\frac{4C_2}{(n-1)^2} \int_0^{C_0\widetilde{K}R^{-\frac{n-1}{2}}\varepsilon}\frac{\mu(\tau)}{\tau}\,\mathrm{d}\tau=+\infty
\end{align*}
 since the modulus of continuity $\mu(\tau)$ satisfies \eqref{Assumption-Thm-Blow-up}.
Our proof of Theorem \ref{Thm-Blow-up-sub} is complete.

\subsection{Proof of Proposition \ref{Prop-Lifespan}}
\hspace{5mm}In order to derive our desired lifespan estimate, let us recall that
\begin{align*}
\ml{H}(\tau)=\int_{\tau}^{\tau_0}\frac{\mu(\sigma)}{\sigma}\,\mathrm{d}\sigma.
\end{align*}
From the previous definition, we see that $\ml{H}$ is decreasing and $\displaystyle{\lim_{\tau\to 0^+}\ml{H}(\tau)=+\infty}$ thanks to the assumption \eqref{Assumption-Thm-Blow-up}. Consequently, $\ml{H}:(0,\tau_0]\to[0,+\infty)$ is invertible.

 From \eqref{Final-Ineq} we derive the following inequality:
\begin{align}\label{Lower-Bound-G}
[\ml{G}(z)]^{-\frac{2}{n-1}} \leqslant (C_0\varepsilon)^{-\frac{2}{n-1}}-\widetilde{K}_1\left(\ml{H}(\widetilde{K}_2z^{-\frac{n-1}{2}}\varepsilon)-\ml{H}(c_2\varepsilon)\right)
\end{align}
for $z\geqslant R$, where $c_2:=C_0\widetilde{K}R^{-\frac{n-1}{2}}$,  $\widetilde{K}_1:= \frac{4C_2}{(n-1)^{2}}$ and $\widetilde{K}_2=C_0\widetilde{K}$. Consequently, for any $z\geqslant R$ such that
\begin{align}\label{impossible z}
(C_0\varepsilon)^{-\frac{2}{n-1}}-\widetilde{K}_1\left(\ml{H}(\widetilde{K}_2z^{-\frac{n-1}{2}}\varepsilon)-\ml{H}(c_2\varepsilon)\right)< 0
\end{align} the inequality in \eqref{Lower-Bound-G} is violated (we recall that according to \eqref{G lb epsilon}, $\ml{G}(z)$ is a positive function), therefore, $\ml{G}(z)$, and in turn $\mathcal{U}(R+z,z)$, does not exist. We remind that we are working on the characteristic line $t=R+z$ for $z\geqslant R$, so, $z\leqslant t= z+R\leqslant 2z$. Since the left-hand side of \eqref{impossible z} is decreasing with respect to $z$, we obtain that $\mathcal{U}(R+z,z)$ cannot be defined whenever $t\geqslant 2R$ fulfills
\begin{align}\label{impossible t}
(C_0\varepsilon)^{-\frac{2}{n-1}}-\widetilde{K}_1\left[\ml{H}\left(\widetilde{K}_2\left(\frac{t}{2}\right)^{-\frac{n-1}{2}}\varepsilon\right)-\ml{H}(c_2\varepsilon)\right]< 0.
\end{align} Then, the lifespan $ T_\varepsilon$ must satisfy
\begin{align} \label{upper bound lifespan implicit}
(C_0\varepsilon)^{-\frac{2}{n-1}}-\widetilde{K}_1\left[\ml{H}\left(\widetilde{K}_3T_\varepsilon^{-\frac{n-1}{2}}\varepsilon\right)-\ml{H}(c_2\varepsilon)\right]\geqslant 0,
\end{align} where $\widetilde{K}_3:=2^{\frac{n-1}{2}}\widetilde{K}_2$. Rearranging \eqref{upper bound lifespan implicit}, we conclude 
\begin{align*}
&(C_0\varepsilon)^{-\frac{2}{n-1}}-\widetilde{K}_1\left[\ml{H}\left(\widetilde{K}_3T_\varepsilon^{-\frac{n-1}{2}}\varepsilon\right)-\ml{H}(c_2\varepsilon)\right]\geqslant 0 \\
&\qquad  \Leftrightarrow \  \widetilde{K}_1^{-1}(C_0\varepsilon)^{-\frac{2}{n-1}} +\ml{H}(c_2\varepsilon)\geqslant \ml{H}(\widetilde{K}_3T_\varepsilon^{-\frac{n-1}{2}}\varepsilon) \\
&\qquad  \Leftrightarrow \  \ml{H}^{-1}\left(\widetilde{K}_1^{-1}(C_0\varepsilon)^{-\frac{2}{n-1}} +\ml{H}(c_2\varepsilon)\right) \leqslant \widetilde{K}_3T_\varepsilon^{-\frac{n-1}{2}}\varepsilon \\
&\qquad  \Leftrightarrow \ (\widetilde{K}_3\varepsilon)^{\frac{2}{n-1}}\left[\ml{H}^{-1}\left(\widetilde{K}_1^{-1}(C_0\varepsilon)^{-\frac{2}{n-1}} +\ml{H}(c_2\varepsilon)\right) \right]^{-\frac{2}{n-1}}\geqslant T_\varepsilon,
\end{align*} which is \eqref{Upper-Lifespan} with $c_1:=\widetilde{K}_1^{-1}C_0^{-\frac{2}{n-1}}$ and $C:=\widetilde{K}_3^{\frac{2}{n-1}}$.
We stress that  \eqref{Upper-Lifespan} is valid only for $\varepsilon\in (0,\varepsilon_0]$ for some $\varepsilon_0=\varepsilon_0(n,R,u_1)>0$. Indeed, in \eqref{impossible t} we have to assume $t\geqslant 2R$ and this condition is compatible with \eqref{Upper-Lifespan} only if we work with $\varepsilon\in (0,\varepsilon_0]$. 

\section{Existence results in the radially symmetric and 3-d case}\label{Section-GESDS-3D}\setcounter{equation}{0}

\hspace{5mm}Before proving the well-posedness results in the radially symmetric and 3-dimensional case, we provide some preliminary results and we explain our approach. We can reformulate the semilinear Cauchy problem \eqref{Wave_Modulus_ut_rad 3d} equivalently as an integral equation.
We begin with the following remark: if $u=u(t,r)$ is a classical solution to \eqref{Wave_Modulus_ut_rad 3d} such that $u\in\mathcal{C}^1([0,+\infty)^2)$, $ru(t,r)\in\mathcal{C}^2([0,+\infty)^2)$ and $u_r(t,0)=0$ for any $t\geqslant 0$, then, we can prolong $u$ to an even function with respect to the variable $r$ such that $u\in\mathcal{C}^1([0,+\infty)\times\mathbb{R})$ and $ru(t,r)\in\mathcal{C}^2([0,+\infty)\times \mathbb{R})$. Moreover, the function $ru$ solves the Cauchy problem
\begin{align}\label{Wave_Modulus_ut_rad 3d not sing}
\begin{cases}
(ru)_{tt}-(ru)_{rr}=r|u_t|^{2}\mu(|u_t|), & r\in\mathbb{R},\ t>0,\\
ru(0,r)=\varepsilon rf_0(r), \ ru_t(0,r)=\varepsilon  rf_1(r) ,&r\in\mathbb{R},
\end{cases}
\end{align} where we keep denoting by $f_0,f_1$ the even extension of $f_0,f_1$ on $\mathbb{R}$. Notice that $rf_0(r)\in\mathcal{C}^{2}(\mathbb{R})$ (thanks to $f_0'(0)=0$) and $rf_1(r)\in\mathcal{C}^{1}(\mathbb{R})$. By D'Alembert's formula, we can represent the solution to \eqref{Wave_Modulus_ut_rad 3d not sing} as
\begin{align*}
u(t,r) = \varepsilon u_{\mathrm{lin}}(t,r)+ u_{\mathrm{nlin}}(t,r),
\end{align*} where 
\begin{align*}
u_{\mathrm{lin}}(t,r) & :=\frac{1}{2r} [(r+t)f_0(r+t)+(r-t)f_0(r-t)]+\frac{1}{2r}  \int_{r-t}^{r+t}\varrho\, f_1(\varrho)\,\mathrm{d}\varrho, \\
u_{\mathrm{nlin}}(t,r) & :=\frac{1}{2r} \int_0^t \int_{r-t+s}^{r+t-s}\varrho\, |u_t(s,\varrho)|^2 \mu(|u_t(s,\varrho)|)\,\mathrm{d}\varrho\, \mathrm{d}s.
\end{align*}  
By introducing the function 
\begin{align}\label{def H function u lin}
H(\lambda):= \lambda f_0(\lambda) +\int_0^\lambda  \varrho f_1(\varrho) \,\mathrm{d}\varrho,
\end{align} we can express $u_{\mathrm{lin}}$ as follows:
\begin{align}\label{representation u lin}
u_{\mathrm{lin}}(t,r) = \frac{1}{2r} [H(r+t)-H(t-r)],
\end{align} where we used that $rf_1(r)$ is an odd function in \eqref{def H function u lin}.
Moreover, we set
\begin{align}\label{def operator I}
I[u](t,r):=\frac{1}{2}  \int_0^t \int_{r-t+s}^{r+t-s}\varrho\, |u_t(s,\varrho)|^2 \mu(|u_t(s,\varrho)|)\,\mathrm{d}\varrho\, \mathrm{d}s
\end{align}  for any $t\geqslant 0$ and any $r\in\mathbb{R}$.

 Consequently, a function $u$ is a  classical solution to \eqref{Wave_Modulus_ut_rad 3d not sing} if and only if $u$ solves the integral equation $u= \varepsilon u_{\mathrm{lin}}+r^{-1} I[u]$, or, equivalently, if $u$ is a fixed point of the nonlinear integral operator $N$ defined by
 \begin{align*}
 N[u](t,r):= \varepsilon u_{\mathrm{lin}}(t,r)+ \frac{1}{r} I[u](t,r)
 \end{align*}
for any $t\geqslant 0$ and $r\in\mathbb{R}$.  Let us introduce now the weighted $L^\infty_{t,r}$ spaces on which we will look for a fixed point of the operator $N$.

\begin{defn}\label{Definition X(R)} Let $R>0$. We consider the weighted space $\mathrm{X}(R)$ such that $u\in \mathrm{X}(R)$ provided that the following conditions hold:
\begin{align}
&\bullet \ u\in \mathcal{C}^1\left([0,+\infty)\times\mathbb{R}\right) \ \mbox{and}\, \ ru(t,r)\in \mathcal{C}^2\left([0,+\infty)\times\mathbb{R}\right),  \label{X(R) regularity} \\
&\bullet \ \forall t\geqslant 0,\  \forall r\in\mathbb{R}: u(t,-r)=u(t,r), \label{X(R) u even in r} \\
&\bullet \ \forall t\geqslant 0, \ \forall r\in\mathbb{R} \ \mbox{such that} \ |r|\geqslant t+R: u(t,r)=0, \label{X(R) supp u} \\
&\bullet \ \|u\|_{\mathrm{X}(R)}:= \sum_{|\alpha|=0,1}\| \phi(t,r) \partial^\alpha u(t,r)\|_{L^{\infty}([0,+\infty)\times \mathbb{R})}+ \sum_{|\alpha|=1,2}\| \phi(t,r) \partial^\alpha(r u(t,r))\|_{L^{\infty}([0,+\infty)\times\mathbb{R})} < +\infty \notag, 
\end{align} where the weight function is given by 
\begin{align}\label{def weight phi}
\phi(t,r):= t-|r|+2R.
\end{align}
\end{defn}

\begin{remark} In the previous definition, \eqref{X(R) supp u} states the property of finite speed propagation for the solution; indeed, we recall that $\forall r\in\mathbb{R}$ such that $|r|\geqslant R$: $f_0(r)=f_1(r)=0$. The choice of the weight function in \eqref{def weight phi} follows the critical case in \cite{Sideris=1983}.
\end{remark}

In Definition \ref{Definition X(R)}, we introduced the space for the global (in time) existence result in Theorem \ref{Thm GESDS rad 3d}, in the next definition, we introduce the local version that we employ in the proof of  Proposition \ref{Prop lower bound lifespan}.

\begin{defn} Let $R>0$ and $T>0$. We consider the weighted space $\mathrm{X}_T(R)$ such that $u\in \mathrm{X}_T(R)$ provided that the following conditions hold:
\begin{align*}
&\bullet \ u\in \mathcal{C}^1\left([0,T]\times\mathbb{R}\right) \ \mbox{and}\, \ ru(t,r)\in \mathcal{C}^2\left([0,T]\times\mathbb{R}\right),   \\
&\bullet \ \forall t\in [0,T],\  \forall r\in\mathbb{R}: u(t,-r)=u(t,r), \\
&\bullet \ \forall t \in [0,T],\ \forall r\in\mathbb{R} \ \mbox{such that} \ |r|\geqslant t+R: u(t,r)=0,  \\
&\bullet \ \|u\|_{\mathrm{X}_T(R)}:= \sum_{|\alpha|=0,1}\| \phi(t,r) \partial^\alpha u(t,r)\|_{L^{\infty}([0,T]\times \mathbb{R})}+ \sum_{|\alpha|=1,2}\| \phi(t,r) \partial^\alpha(r u(t,r))\|_{L^{\infty}([0,T]\times\mathbb{R})} < +\infty, 
\end{align*} where $\phi(t,r)$ is defined in
\eqref{def weight phi}.
\end{defn}

By using \eqref{representation u lin} and the support condition for $f_0,f_1$, one can prove that $u_{\mathrm{lin}} \in X(R)$ and obtain the estimate
\begin{align*}
\|u_{\mathrm{lin}}\|_{X(R)} & \lesssim \sum_{|\alpha|=0,1}\| \partial^\alpha u_{\mathrm{lin}}(t,r)\|_{L^{\infty}([0,+\infty)\times \mathbb{R})}+ \sum_{|\alpha|=1,2}\| \partial^\alpha(r u_{\mathrm{lin}}(t,r))\|_{L^{\infty}([0,+\infty)\times\mathbb{R})} \\ & \lesssim \| (H', H'')\|_{L^{\infty}(\mathbb{R})} 
\end{align*} (see \cite[Equations (10)-(14)]{Sideris=1983} for further details), which implies that
\begin{align}
\|u_{\mathrm{lin}}\|_{X(R)}  & \leqslant  \widetilde{c}_0\|(f_0,f'_0,f''_0,f_1,f'_1)\|_{L^{\infty}([0,+\infty))}, \label{norm u lin}
\end{align} where $\widetilde{c}_0=\widetilde{c}_0(R)>0$.

In the next subsections, we are going to prove that $N$ has a uniquely determined fixed point in $\mathrm{X}(R)$ and $\mathrm{X}_T(R)$ under the assumptions of Theorem \ref{Thm GESDS rad 3d} and Proposition \ref{Prop lower bound lifespan}, respectively. Therefore, to achieve this goal by applying the Banach contraction principle, we have to show that following two fundamental inequalities:
\begin{align}
\| N[u]\|_{\mathrm{X}_T(R)} & \leqslant \widetilde{c}_0(R)  \|(f_0,f'_0,f''_0,f_1,f'_1)\|_{L^{\infty}([0,+\infty))} \, \varepsilon+ \widetilde{c}(R,\mu,\delta,T) \|u\|^2_{\mathrm{X}_T(R)}, \label{1st ineq contraction} \\
\| N[u] -N[v]\|_{\mathrm{X}_T(R)} &  \leqslant \widetilde{c}(R,\mu,\delta,T) \|u-v\|_{\mathrm{X}_T(R)} \left(\|u\|_{\mathrm{X}_T(R)}+\|v\|_{\mathrm{X}_T(R)}\right), \label{2nd ineq contraction}
\end{align} 
hold for any $u,v$ such that $\|u\|_{\mathrm{X}_T(R)},\|v\|_{\mathrm{X}_T(R)}\leqslant \delta$, where $\delta$ is a positive constant to be fixed later. More precisely, in Subsection \ref{Subsection GESDS 3d rad} we use the Dini condition in \eqref{Dini condition} for $\mu$ to prove that $\widetilde{c}(R,\mu,\delta,T)$ is uniformly bounded with respect to the time variable, 
while in Subsection \ref{Subsection lower bound lifespan 3d rad} we express $\widetilde{c}(R,\mu,\delta,T)$ through the function $\mathcal{H}$ defined in \eqref{def H function lifespan} to prove that a local solution to \eqref{Wave_Modulus_ut_rad 3d not sing} exists at least for $T$ satisfying \eqref{lower bound lifespan 3d rad-T}.

\subsection{Proof of Theorem \ref{Thm GESDS rad 3d}}
\label{Subsection GESDS 3d rad}
\hspace{5mm}We begin by showing that for any $u\in\mathrm{X}(R)$ we have $r^{-1} I[u](t,r)\in \mathrm{X}(R)$.\\

 Since $u_t(t,\cdot)$ is an even function, we have $I[u](t,-r)=-I[u](t,r)$ for any $t\geqslant 0$ and $r\in\mathbb{R}$. Hence, $r^{-1} I[u](t,r)$  satisfies \eqref{X(R) u even in r} and
\begin{align} \label{I[u](t,0)=0}
I[u](t,0)=0 \qquad \mbox{for any} \ t\geqslant 0.
\end{align}

Furthermore, for any $t\geqslant 0$ and $r\in\mathbb{R}$ such that $|r|\geqslant t+R$ we get $I[u](t,r)=0$, since the integrand function in \eqref{def operator I} is identically 0 in this region, so $r^{-1}I[u](t,r)$ satisfies \eqref{X(R) supp u} as well.\\

 Let us prove now that $I[u]\in\mathcal{C}^2([0,+\infty)\times\mathbb{R})$. Since $\mu$ and $u_t$ are continuous, from \eqref{def operator I} we have that $I[u]$ is continuous. Let us calculate the derivatives $\partial_r I[u]$, $\partial_t I[u]$, $\partial_r^2 I[u]$, $\partial_t \partial_r I[u]$, $\partial_t^2 I[u]$ and show that these are continuous. Denoting $G(s,\varrho):=|u_t(s,\varrho)|^2\mu(|u_t(s,\varrho)|)$ for the sake of brevity, by a direct computation we find
\begin{align*}
\partial_r I[u](t,r)=\frac{1}{2}\int_0^t (r+t-s)G(s,r+t-s)\, \mathrm{d}s-\frac{1}{2}\int_0^t (r-t+s)G(s,r-t+s)\, \mathrm{d}s.
\end{align*} From the support condition \eqref{X(R) supp u}, it follows that
\begin{align*}
& s\leqslant \tfrac{1}{2}(t+r-R) \quad \Rightarrow \quad r+t-s\geqslant s+R \quad\ \ \quad \Rightarrow \quad  G(s,r+t-s)=0, \\
& s\leqslant \tfrac{1}{2}(t-r-R) \quad \Rightarrow \quad r-t+s\leqslant -(s+R) \quad \Rightarrow \quad  G(s,r-t+s)=0.
\end{align*} Therefore,
\begin{align}
\partial_r I[u](t,r)=\frac{1}{2}\int_{t_1}^t \varrho \, G(s,\varrho)\big|_{\varrho=r+t-s} \mathrm{d}s-\frac{1}{2}\int_{t_2}^t \varrho \, G(s,\varrho)\big|_{\varrho=r-t+s} \mathrm{d}s,\label{der r I[u]}
\end{align} where
\begin{align*}
t_1:= \max\{0,\tfrac{1}{2}(t+r-R)\}, \quad t_2:= \max\{0,\tfrac{1}{2}(t-r-R)\}.
\end{align*} Analogously,
\begin{align}
\partial_t I[u](t,r)=\frac{1}{2}\int_{t_1}^t \varrho \, G(s,\varrho)\big|_{\varrho=r+t-s} \mathrm{d}s+\frac{1}{2}\int_{t_2}^t \varrho \, G(s,\varrho)\big|_{\varrho=r-t+s} \mathrm{d}s.\label{der t I[u]}
\end{align} Since $\varrho\, G(s,\varrho)$ is a continuous function, from \eqref{der r I[u]} and \eqref{der t I[u]} we have that $\partial_r I[u],\partial_t I[u]$ are continuous too.
Similar computations lead to the representations
\begin{align}
\partial_r^2 I[u](t,r) & =\frac{1}{2}\int_{t_1}^t \frac{\partial}{\partial \varrho}\big(\varrho \, G(s,\varrho)\big)\big|_{\varrho=r+t-s} \mathrm{d}s-\frac{1}{2}\int_{t_2}^t \frac{\partial}{\partial \varrho}\big(\varrho \, G(s,\varrho)\big)\big|_{\varrho=r-t+s} \mathrm{d}s,\label{der rr I[u]} \\
\partial_t \partial_r I[u](t,r) & =\frac{1}{2}\int_{t_1}^t \frac{\partial}{\partial \varrho}\big(\varrho \, G(s,\varrho)\big)\big|_{\varrho=r+t-s} \mathrm{d}s+\frac{1}{2}\int_{t_2}^t \frac{\partial}{\partial \varrho}\big(\varrho \, G(s,\varrho)\big)\big|_{\varrho=r-t+s} \mathrm{d}s,\label{der tr I[u]} \\
\partial_t^2 I[u](t,r) &=  \partial_r^2 I[u](t,r) +r G(t,r). \label{der tt I[u]}
\end{align}
We denote by $F(\tau):=|\tau|^2\mu(|\tau|)$, so that $G(s,\varrho)=F(u_t(s,\varrho))$. Thanks to the regularity assumption on $\mu$ in $(A_3)$ of Theorem \ref{Thm GESDS rad 3d} we have that $F\in\mathcal{C}^2(\mathbb{R})$ and 
\begin{equation}\label{der F nonlinearity}
\begin{split} 
F'(\tau) &= 2\tau\mu(|\tau|)+\tau|\tau|\mu'(|\tau|), \\
F''(\tau) &= 2\mu(|\tau|)+4|\tau|\mu'(|\tau|)+|\tau|^2\mu''(|\tau|).
\end{split}
\end{equation} 
Then, 
\begin{align*}
\frac{\partial}{\partial\varrho}\big(\varrho\, G(s,\varrho)\big)=G(s,\varrho)+ F'(u_t(s,\varrho)) \varrho\, u_{tr}(s,\varrho)
\end{align*}
 is a continuous function (notice that $r u_{tr}=(ru)_{tr}-u_t$ is continuous thanks to \eqref{X(R) regularity}). Consequently, from \eqref{der rr I[u]},   \eqref{der tr I[u]},  \eqref{der tt I[u]} we have that $\partial_r^2 I[u],\partial_t \partial_r I[u],\partial_t^2 I[u]$ are continuous functions.
Since we want to show that $r^{-1}I[u](t,r)$ satisfies \eqref{X(R) regularity}, it remains to prove that $r^{-1}I[u](t,r)\in \mathcal{C}^1([0,+\infty)\times \mathbb{R})$. By using \eqref{I[u](t,0)=0}, we obtain
\begin{align}
\frac{1}{r} I[u](t,r) & = \frac{1}{r} \int_0^r \partial_r I[u](t,\rho)\,\mathrm{d}\rho \label{1/r I[u]}, \\
\frac{\partial}{\partial t}\left(\frac{1}{r}  I[u](t,r)\right) & = \frac{1}{r} \int_0^r \partial_t \partial_r I[u](t,\rho)\,\mathrm{d}\rho \label{1/r I[u] der t}, \\
\frac{\partial}{\partial r}\left(\frac{1}{r}  I[u](t,r)\right) & = \frac{1}{r^2} \int_0^r \rho \, \partial_r^2 I[u](t,\rho)\,\mathrm{d}\rho \label{1/r I[u] der r},
\end{align} where for the $r$-derivative we integrate by parts after differentiating \eqref{1/r I[u]} with respect to $r$. Notice that $\partial_r(r^{-1}I[u](t,r))\big|_{r=0}=0$. Since we proved that $I[u]\in\mathcal{C}^2([0,+\infty)\times\mathbb{R})$, from \eqref{1/r I[u]}, \eqref{1/r I[u] der t} and  \eqref{1/r I[u] der r} we conclude that $r^{-1}I[u](t,r)\in \mathcal{C}^1([0,+\infty)\times \mathbb{R})$.\\

 Finally, in order to show that $r^{-1}I[u](t,r)\in \mathrm{X}(R)$ we have to estimate $\|r^{-1} I[u](t,r)\|_{\mathrm{X}(R)}$.
By the elementary relations
\begin{align*}
u_r&=-\tfrac{1}{r} u+ \tfrac{1}{r} \partial_r(ru),  &&u_{rr}=-\tfrac{2}{r} u_r+ \tfrac{1}{r} \partial_r^2(ru), \\ 
u_t& =\tfrac{1}{r}\partial_t(ru), &&u_{tt}=\tfrac{1}{r}\partial_t^2(ru),  &&u_{tr}=-\tfrac{1}{r} u_t+\tfrac{1}{r}\partial_t\partial_r(ru),
\end{align*} and the definition of the norm $\|\cdot\|_{\mathrm{X}(R)}$, we derive the upper bounds
\begin{align}
|\partial^\alpha u(t,r)| & \leqslant \phi(t,r)^{-1} \|u\|_{\mathrm{X}(R)} && \mbox{for} \ |\alpha|=0,1, \label{upper partial u |alpha|=0,1}\\
|\partial^\alpha u(t,r)| & \leqslant \widetilde{C}|r|^{-1} \phi(t,r)^{-1} \|u\|_{\mathrm{X}(R)} && \mbox{for} \ |\alpha|=1,2, \label{upper partial u |alpha|=1,2} \\
|\partial^\alpha u(t,r)| & \leqslant \widetilde{C} \vartheta(r) \phi(t,r)^{-1} \|u\|_{\mathrm{X}(R)} && \mbox{for} \ |\alpha|=1, \label{upper partial u |alpha|=1}
\end{align} for any $t\geqslant 0$ and $r\in\mathbb{R}$, where $\vartheta(r):=\min\left\{1,\frac{R}{|r|}\right\}$ and $\widetilde{C}=\widetilde{C}(R)\geqslant1$.
Therefore, combining \eqref{mu derivatives growth cond}, \eqref{der F nonlinearity} with \eqref{upper partial u |alpha|=0,1}, \eqref{upper partial u |alpha|=1,2}, \eqref{upper partial u |alpha|=1} we have
\begin{align}
|\varrho\, G(s,\varrho)| & = |\varrho|\, |u_t(s,\varrho)|^2 \mu(|u_t(s,\varrho)|)\notag\\
& \lesssim \vartheta(\varrho)\phi(s,\varrho)^{-2} \mu\left(\widetilde{C}\|u\|_{\mathrm{X}(R)}\vartheta(\varrho)\phi(s,\varrho)^{-1}\right) \|u\|^2_{\mathrm{X}(R)}, \label{upper bound rhoG(s,rho)} \\
| \partial_{\varrho}(\varrho \, G(s,\varrho))|  & = |G(s,\varrho)+\varrho F'(u_t(s,\varrho)) u_{tr}(s,\varrho)| \notag \\
 & \lesssim |u_t(s,\varrho)|^2 \mu(|u_t(s,\varrho)|)+|\varrho \, u_{tr}(s,\varrho)|\,|u_t(s,\varrho)|\,  \mu(|u_t(s,\varrho)|) \notag\\ 
& \lesssim \vartheta(\varrho)\phi(s,\varrho)^{-2} \mu\left(\widetilde{C}\|u\|_{\mathrm{X}(R)}\vartheta(\varrho)\phi(s,\varrho)^{-1}\right) \|u\|^2_{\mathrm{X}(R)}.\label{upper bound rhoG(s,rho) der rho}
\end{align} By \eqref{der r I[u]}, \eqref{der t I[u]} and \eqref{upper bound rhoG(s,rho)} we obtain
\begin{align*}
|\partial_r I[u](t,r)|,|\partial_t I[u](t,r)| \lesssim \big(J_1(t,r)+J_2(t,r)\big) \|u\|^2_{\mathrm{X}(R)},
\end{align*} where
\begin{align}
J_1(t,r) & := \int_{t_1}^t  \vartheta(\varrho)\phi(s,\varrho)^{-2} \mu\left(\widetilde{C}\|u\|_{\mathrm{X}(R)}\vartheta(\varrho)\phi(s,\varrho)^{-1}\right) \big|_{\varrho=r+t-s} \mathrm{d}s, \label{def J1(t,r)}\\
J_2(t,r) & := \int_{t_2}^t  \vartheta(\varrho)\phi(s,\varrho)^{-2} \mu\left(\widetilde{C}\|u\|_{\mathrm{X}(R)}\vartheta(\varrho)\phi(s,\varrho)^{-1}\right) \big|_{\varrho=r-t+s} \mathrm{d}s. \label{def J2(t,r)}
\end{align} Similarly, combining \eqref{der rr I[u]}, \eqref{der tr I[u]}, \eqref{der tt I[u]}, \eqref{upper bound rhoG(s,rho)} and \eqref{upper bound rhoG(s,rho) der rho}, we get
\begin{align*}
|\partial_r^2 I[u](t,r)|,|\partial_t\partial_r I[u](t,r)| & \lesssim \big(J_1(t,r)+J_2(t,r)\big) \|u\|^2_{\mathrm{X}(R)}, \\
|\partial_t^2 I[u](t,r)| &\lesssim \left(J_1(t,r)+J_2(t,r)+\phi(t,r)^{-2}\mu\left(\|u\|_{\mathrm{X}(R)}\phi(t,r)^{-1}\right)\right) \|u\|^2_{\mathrm{X}(R)}.
\end{align*}

\begin{lemma} \label{Lemma J1,J2} For any $(t,r)\in[0,+\infty)\times\mathbb{R}$ such that $|r|\leqslant t+R$ the following estimates hold:
\begin{align}\label{estimare J1,J2 lemma}
J_h(t,r) \lesssim \left(\mu\left(\widetilde{C}R^{-1}\|u\|_{\mathrm{X}(R)}\right) +\int_{\widetilde{C}(t+R)^{-1}\|u\|_{\mathrm{X}(R)}}^{\widetilde{C}R^{-1}\|u\|_{\mathrm{X}(R)}} \frac{\mu(\tau)}{\tau}\, \mathrm{d}\tau\right) \phi(t,r)^{-1} \qquad \mbox{for} \ h=1,2,
\end{align} 
where the unexpressed multiplicative constant depends only on $R$.
\end{lemma}

\noindent We postpone the proof of Lemma \ref{Lemma J1,J2} to the end of this subsection. We employ now the Dini condition for $\mu$: from \eqref{Dini condition} it follows that the factor between the round brackets in the right-hand side of \eqref{estimare J1,J2 lemma} is finite and can be estimated from above by 
\begin{align*}
M(R,\mu,\|u\|_{\mathrm{X}(R)}):=\mu\left(\widetilde{C}R^{-1}\|u\|_{\mathrm{X}(R)}\right) +\int_0^{\widetilde{C}R^{-1}\|u\|_{\mathrm{X}(R)}} \frac{\mu(\tau)}{\tau}\, \mathrm{d}\tau.
\end{align*} 
Summarizing,  since $\phi(t,r)\geqslant R$, we proved that 
\begin{align} \label{upper bound partial I[u] |alpha|=1,2}
|\partial^\alpha I[u](t,r)| \lesssim M(R,\mu,\|u\|_{\mathrm{X}(R)}) \|u\|_{\mathrm{X}(R)}^2 \phi(t,r)^{-1}
\end{align} for any $\alpha\in\mathbb{N}^2$ with $|\alpha|=1,2$ and any $(t,r)\in[0,+\infty)\times \mathbb{R}$ such that $|r|\leqslant t+R$, where $\partial=(\partial_t,\partial_r)$. Therefore, from \eqref{upper bound partial I[u] |alpha|=1,2} and \eqref{1/r I[u]}
\begin{align} 
\left|\frac{1}{r}I[u](t,r)\right|  \leqslant \frac{1}{|r|} \int_0^{|r|} |\partial_r I[u](t,\rho)|\,\mathrm{d}\rho &\lesssim M(R,\mu,\|u\|_{\mathrm{X}(R)}) \|u\|_{\mathrm{X}(R)}^2  \frac{1}{|r|} \int_0^{|r|} \phi(t,\rho)^{-1}\, \mathrm{d}\rho \notag \\
& \lesssim M(R,\mu,\|u\|_{\mathrm{X}(R)}) \|u\|_{\mathrm{X}(R)}^2 \phi(t,r)^{-1},  \label{upper bound partial 1/r I[u]}
\end{align}
where we used that $\phi(t,r)\leqslant\phi(t,\rho)$ for any $\rho\in\mathbb{R}$ such that $|\rho|\leqslant |r|$. In a similar way, by \eqref{1/r I[u] der t}, \eqref{1/r I[u] der r} we can prove that 
\begin{align}
\left|\frac{\partial}{\partial t}\left(\frac{1}{r}I[u](t,r)\right)\right|,\left|\frac{\partial}{\partial r}\left(\frac{1}{r}I[u](t,r)\right)\right| \lesssim M(R,\mu,\|u\|_{\mathrm{X}(R)}) \|u\|_{\mathrm{X}(R)}^2 \phi(t,r)^{-1}. \label{upper bound partial 1/r I[u] |alpha|=1}
\end{align}
Combining \eqref{upper bound partial I[u] |alpha|=1,2}, \eqref{upper bound partial 1/r I[u]} and \eqref{upper bound partial 1/r I[u] |alpha|=1} we proved that 
\begin{align}
\| r^{-1} I[u](t,r)\|_{\mathrm{X}(R)}\leqslant \widetilde{C}_1 \left(\mu\left(\widetilde{C}R^{-1}\|u\|_{\mathrm{X}(R)}\right) +\int_0^{\widetilde{C}R^{-1}\|u\|_{\mathrm{X}(R)}} \frac{\mu(\tau)}{\tau}\, \mathrm{d}\tau\right) \|u\|_{\mathrm{X}(R)}^2, \label{norm 1/r I[u]}
\end{align} where the multiplicative constant $\widetilde{C}_1>0$ depends only on $R$. From \eqref{norm u lin} and \eqref{norm 1/r I[u]}, we obtain the desired estimate \eqref{1st ineq contraction}.

Next we derive \eqref{2nd ineq contraction}.
From \eqref{der r I[u]} and \eqref{der t I[u]} we have
\begin{align*}
\partial_r \big\{I[u](t,r)-I[v](t,r)\big\} &= \frac{1}{2}\int_{t_1}^t \varrho\, \big( G(s,\varrho)-\widetilde{G}(s,\varrho)\big)\big|_{\varrho=r+t-s} \mathrm{d}s-\frac{1}{2}\int_{t_2}^t \varrho\, \big( G(s,\varrho)-\widetilde{G}(s,\varrho)\big)\big|_{\varrho=r-t+s} \mathrm{d}s, \\
\partial_t \big\{I[u](t,r)-I[v](t,r)\big\} &= \frac{1}{2}\int_{t_1}^t \varrho\, \big( G(s,\varrho)-\widetilde{G}(s,\varrho)\big)\big|_{\varrho=r+t-s} \mathrm{d}s+\frac{1}{2}\int_{t_2}^t \varrho\, \big( G(s,\varrho)-\widetilde{G}(s,\varrho)\big)\big|_{\varrho=r-t+s} \mathrm{d}s,
\end{align*} where $G(s,\varrho)=F(u_t(s,\varrho))$ and $\widetilde{G}(s,\varrho)=F(v_t(s,\varrho))$.

 By \eqref{der F nonlinearity} and \eqref{mu derivatives growth cond}, it follows that $|F'(\tau)|\lesssim |\tau|\mu(|\tau|)$. Applying the mean value theorem, the monotonicity of $\mu$, \eqref{upper partial u |alpha|=1,2} and \eqref{upper partial u |alpha|=1}, we arrive at
\begin{align}
|\varrho|\, \big|G(s,\varrho) & -\widetilde{G}(s,\varrho)\big| = |\varrho| \,\big|F'\big(\theta_1 u_t(s,\varrho)+(1-\theta_1) v_t(s,\varrho)\big)\big|\, |u_t(s,\varrho)-v_t(s,\varrho)| \notag\\
 & \lesssim |\varrho| \left(|u_t(s,\varrho)|+|v_t(s,\varrho)|\right) \mu\left(|u_t(s,\varrho)|+|v_t(s,\varrho)|\right)|u_t(s,\varrho)-v_t(s,\varrho)| \notag \\
 & \lesssim \vartheta(\varrho)\phi(s,\varrho)^{-2} \mu\big(\widetilde{C}(\|u\|_{\mathrm{X}(R)}+\|v\|_{\mathrm{X}(R)})\vartheta(\varrho)\phi(s,\varrho)^{-1}\big) \|u-v\|_{\mathrm{X}(R)} \left(\|u\|_{\mathrm{X}(R)}+\|v\|_{\mathrm{X}(R)}\right), \label{estimate rho G-Gtilde}
\end{align} where in the first step we took some $\theta_1\in(0,1)$ depending on $s,\varrho,u_t,v_t$. Therefore,
\begin{align*}
 &\big|\partial_r \big\{I[u](t,r)-I[v](t,r)\big\}\big| \lesssim \big(J_3(t,r)+J_4(t,r)\big)\|u-v\|_{\mathrm{X}(R)} \left(\|u\|_{\mathrm{X}(R)}+\|v\|_{\mathrm{X}(R)}\right),\\
 & \big|\partial_t \big\{I[u](t,r)-I[v](t,r)\big\}\big|  \lesssim \big(J_3(t,r)+J_4(t,r)\big)\|u-v\|_{\mathrm{X}(R)} \left(\|u\|_{\mathrm{X}(R)}+\|v\|_{\mathrm{X}(R)}\right),
\end{align*} where
\begin{align*}
J_3(t,r) & := \int_{t_1}^t  \vartheta(\varrho)\phi(s,\varrho)^{-2} \mu\big(\widetilde{C}(\|u\|_{\mathrm{X}(R)}+\|v\|_{\mathrm{X}(R)})\vartheta(\varrho)\phi(s,\varrho)^{-1}\big) \big|_{\varrho=r+t-s} \mathrm{d}s, \\
J_4(t,r) & := \int_{t_2}^t  \vartheta(\varrho)\phi(s,\varrho)^{-2} \mu\big(\widetilde{C}(\|u\|_{\mathrm{X}(R)}+\|v\|_{\mathrm{X}(R)})\vartheta(\varrho)\phi(s,\varrho)^{-1}\big) \big|_{\varrho=r-t+s} \mathrm{d}s.
\end{align*} Clearly, $J_3(t,r), J_4(t,r)$ can be estimated analogously as the terms $J_1(t,r), J_2(t,r)$ are estimated in Lemma \ref{Lemma J1,J2}: the obvious difference is that one must replace $\|u\|_{\mathrm{X}(R)}$ by $\|u\|_{\mathrm{X}(R)}+\|v\|_{\mathrm{X}(R)}$ in the right-hand side of \eqref{estimare J1,J2 lemma}.

Employing \eqref{der rr I[u]}, \eqref{der tr I[u]} and \eqref{der tt I[u]}, we obtain
\begin{align}
\partial_r^2 \big\{I[u](t,r)-I[v](t,r)\big\} & =\frac{1}{2}\int_{t_1}^t \frac{\partial}{\partial \varrho}\left(\varrho \left( G(s,\varrho)-\widetilde{G}(s,\varrho)\right)\right)\big|_{\varrho=r+t-s} \mathrm{d}s \notag \\ & \quad -\frac{1}{2}\int_{t_2}^t \frac{\partial}{\partial \varrho}\left(\varrho \left( G(s,\varrho)-\widetilde{G}(s,\varrho)\right)\right)\big|_{\varrho=r-t+s} \mathrm{d}s, \label{der rr I[u]- I[v]} \\
\partial_t \partial_r \big\{I[u](t,r)-I[v](t,r)\big\} & =\frac{1}{2}\int_{t_1}^t \frac{\partial}{\partial \varrho}\left(\varrho \left( G(s,\varrho)-\widetilde{G}(s,\varrho)\right)\right)\big|_{\varrho=r+t-s} \mathrm{d}s\notag  \\ & \quad +\frac{1}{2}\int_{t_2}^t \frac{\partial}{\partial \varrho}\left(\varrho \left( G(s,\varrho)-\widetilde{G}(s,\varrho)\right)\right)\big|_{\varrho=r-t+s} \mathrm{d}s, \label{der tr I[u]- I[v]} \\
\partial_t^2 \big\{I[u](t,r)-I[v](t,r)\big\} & = \partial_r^2 \big\{I[u](t,r)-I[v](t,r)\big\} +r \left(G(t,r)-\widetilde{G}(t,r)\right). \label{der tt I[u]- I[v]}
\end{align} Since
\begin{align*}
\partial_\varrho \left(\varrho\, ( G(s,\varrho)-\widetilde{G}(s,\varrho))\right)=\left( G(s,\varrho)-\widetilde{G}(s,\varrho)\right)+\varrho\, \partial_\varrho \left(G(s,\varrho)-\widetilde{G}(s,\varrho)\right)
\end{align*} we estimate the two addends on the right-hand side. For the first term, as in \eqref{estimate rho G-Gtilde}, by using \eqref{upper partial u |alpha|=0,1} and \eqref{upper partial u |alpha|=1} we get
\begin{align}
& \big|G(s,\varrho)  -\widetilde{G}(s,\varrho)\big| \notag \\& \quad \lesssim \left(|u_t(s,\varrho)|+|v_t(s,\varrho)|\right) \mu\left(|u_t(s,\varrho)|+|v_t(s,\varrho)|\right)|u_t(s,\varrho)-v_t(s,\varrho)| \notag\\
 & \quad \lesssim \vartheta(\varrho)\phi(s,\varrho)^{-2} \mu\big(\widetilde{C}(\|u\|_{\mathrm{X}(R)}+\|v\|_{\mathrm{X}(R)})\vartheta(\varrho)\phi(s,\varrho)^{-1}\big) \|u-v\|_{\mathrm{X}(R)} \left(\|u\|_{\mathrm{X}(R)}+\|v\|_{\mathrm{X}(R)}\right). \label{estimate G-Gtilde}
\end{align} For the second term, we rewrite
\begin{align*}
& \frac{\partial}{\partial \varrho}\big( G(s,\varrho)-\widetilde{G}(s,\varrho)\big)   \\& \qquad=\frac{\partial}{\partial \varrho} \int_0^1 \frac{\partial}{\partial \omega} \big[F\big(\omega u_t(s,\varrho)+(1-\omega)v_t(s,\varrho)\big) \big]\mathrm{d}\omega \\
& \qquad =\frac{\partial}{\partial \varrho} \left[\int_0^1 F'\big(\omega u_t(s,\varrho)+(1-\omega)v_t(s,\varrho)\big) \, \mathrm{d}\omega \, \big(u_t(s,\varrho) - v_t(s,\varrho)\big) \right]  \\
& \qquad =\int_0^1 F''\big(\omega u_t(s,\varrho)+(1-\omega)v_t(s,\varrho)\big)\big(\omega u_{tr}(s,\varrho)+(1-\omega)v_{tr}(s,\varrho)\big)  \, \mathrm{d}\omega \,  \big(u_t(s,\varrho) - v_t(s,\varrho)\big) \\ & \qquad \quad + \int_0^1 F'\big(\omega u_t(s,\varrho)+(1-\omega)v_t(s,\varrho)\big) \,\mathrm{d}\omega \, \big(u_{tr}(s,\varrho) - v_{tr}(s,\varrho)\big).
\end{align*} 
Then, we use the monotonicity of $\mu$, the inequalities $|F'(\tau)|\lesssim |\tau|\mu(|\tau|)$, $|F''(\tau)|\lesssim \mu(|\tau|)$, which follow from \eqref{mu derivatives growth cond} and \eqref{der F nonlinearity}, and the estimates in \eqref{upper partial u |alpha|=1,2} and \eqref{upper partial u |alpha|=1}, obtaining
\begin{align}
&\big|\varrho\, \partial_\varrho \big( G(s,\varrho)-\widetilde{G}(s,\varrho)\big) \big| \notag \\ & \quad \lesssim  \mu\left(|u_t(s,\varrho)|+|v_t(s,\varrho)|\right) |\varrho| \left(|u_{tr}(s,\varrho)|+|v_{tr}(s,\varrho)|\right) |u_t(s,\varrho)-v_t(s,\varrho)| \notag \\ & \quad \quad +\left(|u_t(s,\varrho)|+|v_t(s,\varrho)|\right) \mu\left(|u_t(s,\varrho)|+|v_t(s,\varrho)|\right) |\varrho| \, |u_{tr}(s,\varrho)-v_{tr}(s,\varrho)| \notag \\
 & \quad \lesssim \vartheta(\varrho)\phi(s,\varrho)^{-2} \mu\big(\widetilde{C}(\|u\|_{\mathrm{X}(R)}+\|v\|_{\mathrm{X}(R)})\vartheta(\varrho)\phi(s,\varrho)^{-1}\big) \|u-v\|_{\mathrm{X}(R)} \left(\|u\|_{\mathrm{X}(R)}+\|v\|_{\mathrm{X}(R)}\right). \label{estimate rho (G-Gtilde) der rho}
\end{align}
Combining \eqref{der rr I[u]- I[v]}, \eqref{der tr I[u]- I[v]}, \eqref{der tt I[u]- I[v]}, \eqref{estimate rho G-Gtilde}, \eqref{estimate G-Gtilde} and \eqref{estimate rho (G-Gtilde) der rho}, we arrive at
\begin{align*}
\big|\partial_r^2 \big\{I[u](t,r)-I[v](t,r)\big\}\big| &\lesssim \big(J_3(t,r)+J_4(t,r)\big)\|u-v\|_{\mathrm{X}(R)} \left(\|u\|_{\mathrm{X}(R)}+\|v\|_{\mathrm{X}(R)}\right),   \\
 \big|\partial_t\partial_r \big\{I[u](t,r)-I[v](t,r)\big\}\big| &\lesssim \big(J_3(t,r)+J_4(t,r)\big)\|u-v\|_{\mathrm{X}(R)} \left(\|u\|_{\mathrm{X}(R)}+\|v\|_{\mathrm{X}(R)}\right), 
\end{align*}
and
\begin{align*}
  \big|\partial_r^2 \big\{I[u](t,r)-I[v](t,r)\big\}\big| & \lesssim \big(J_3(t,r)+J_4(t,r)\big)\|u-v\|_{\mathrm{X}(R)} \left(\|u\|_{\mathrm{X}(R)}+\|v\|_{\mathrm{X}(R)}\right) \\ +\phi(t,r)^{-2} & \mu\big((\|u\|_{\mathrm{X}(R)}+\|v\|_{\mathrm{X}(R)})\phi(t,r)^{-1}\big) \|u-v\|_{\mathrm{X}(R)} \left(\|u\|_{\mathrm{X}(R)}+\|v\|_{\mathrm{X}(R)}\right).
\end{align*} Summarizing, we proved that 
\begin{align*}
\big|\partial^\alpha  \big\{I[u](t,r)-I[v](t,r)\big\}\big| \lesssim M(R,\mu,\|u\|_{\mathrm{X}(R)}+\|v\|_{\mathrm{X}(R)})\|u-v\|_{\mathrm{X}(R)} \left(\|u\|_{\mathrm{X}(R)}+\|v\|_{\mathrm{X}(R)}\right)\phi(t,r)^{-1}
\end{align*} for any $\alpha\in\mathbb{N}: |\alpha|=1,2$ and any $(t,r)\in [0,+\infty)\times\mathbb{R}$ such that $|r|\leqslant t+R$. With the same argument used to derive \eqref{norm 1/r I[u]}, applying the previous inequalities in the representations analogous to those in \eqref{1/r I[u]}, \eqref{1/r I[u] der t}, \eqref{1/r I[u] der t} for $r^{-1}(I[u](t,r)-I[v](t,r))$, we find that 
\begin{align*}
& \big|\partial^\alpha  \big\{r^{-1}\big(I[u](t,r)-I[v](t,r)\big)\big\}\big| \\ & \qquad\lesssim M(R,\mu,\|u\|_{\mathrm{X}(R)}+\|v\|_{\mathrm{X}(R)})\|u-v\|_{\mathrm{X}(R)} \left(\|u\|_{\mathrm{X}(R)}+\|v\|_{\mathrm{X}(R)}\right) \phi(t,r)^{-1}
\end{align*} for any $\alpha\in\mathbb{N}: |\alpha|=0,1$ and any $(t,r)\in [0,+\infty)\times\mathbb{R}$ such that $|r|\leqslant t+R$. Altogether, we derived the estimate
\begin{align}
\big\| r^{-1} \big(I[u](t,r)-I[v](t,r)\big)\big\|_{\mathrm{X}(R)}   \leqslant \widetilde{C}_2 M(R,\mu,\|u\|_{\mathrm{X}(R)}+\|v\|_{\mathrm{X}(R)}) \|u-v\|_{\mathrm{X}(R)} \left(\|u\|_{\mathrm{X}(R)}+\|v\|_{\mathrm{X}(R)}\right),\label{norm 1/r (I[u]-I[v])}
\end{align} where $\widetilde{C}_2=\widetilde{C}_2(R)>0$.

Finally, we check that \eqref{norm u lin}, \eqref{norm 1/r I[u]} and \eqref{norm 1/r (I[u]-I[v])} imply that the operator $N$ is a contraction on the ball $\mathfrak{B}_\delta(\mathrm{X}(R)):=\{u\in\mathrm{X}(R): \|u\|_{\mathrm{X}(R)}\leqslant \delta\}$ provided that $\delta=\delta(\varepsilon,f_0,f'_0,f''_0,f_1,f'_1,R,\mu)>0$ and $\varepsilon>0$ are sufficiently small. By using the Dini condition and the continuity of $\mu(\tau)$ as $\tau\to 0^+$, we can pick $\delta>0$ such that 
\begin{align*}
2\mu(2\delta \widetilde{C}R^{-1})\leqslant 1, \quad  2\int_{0}^{2\delta \widetilde{C}R^{-1}} \frac{\mu(\tau)}{\tau} \,\mathrm{d}\tau \leqslant 1, \quad 2\widetilde{C}_1 \delta^2\leqslant \delta, \quad  4\widetilde{C}_2\delta\leqslant  1,
\end{align*} and $\varepsilon$ such that $2\widetilde{c}_0 \varepsilon \|(f_0,f'_0,f''_0,f_1,f'_1)\|_{L^\infty([0,+\infty))} =\delta $. According to this setting,  \eqref{norm u lin}, \eqref{norm 1/r I[u]} imply that  $\|  N[u]\|_{\mathrm{X}(R)}\leqslant \delta$ for any $u\in \mathfrak{B}_\delta(\mathrm{X}(R))$ and \eqref{norm 1/r (I[u]-I[v])} implies that  $\|  N[u]-N[v]\|_{\mathrm{X}(R)}\leqslant \frac{1}{2}\|u-v\|_{\mathrm{X}(R)}$ for any $u,v\in \mathfrak{B}_\delta(\mathrm{X}(R))$. By Banach's fixed point theorem it follows the existence of a unique $u\in \mathfrak{B}_\delta(\mathrm{X}(R))$ satisfying the integral equation equivalent to \eqref{Wave_Modulus_ut_rad 3d not sing}.

\begin{proof}[Proof of Lemma \ref{Lemma J1,J2}]

\textbf{Estimate of $J_1(t,r)$:}

 Let us denote $\xi:=r+t$ and $\eta:=t-r$. From \eqref{def J1(t,r)} we have
\begin{align*}
J_1(t,r) & =\int_{\max\{0,\frac{1}{2}(r+t-R)\}}^t \vartheta(r+t-s)\phi(s,r+t-s)^{-2}\mu\left(\widetilde{C}\|u\|_{\mathrm{X}(R)}\vartheta(r+t-s)\phi(s,r+t-s)^{-1}\right)\mathrm{d}s \\
& =\frac{1}{2}\int_{\max\{-\xi,-R\}}^\eta \vartheta\left(	\tfrac{\xi-\sigma}{2}\right)\phi\left(\tfrac{\sigma+\xi}{2},\tfrac{\xi-\sigma}{2}\right)^{-2}\mu\left(\widetilde{C}\|u\|_{\mathrm{X}(R)}\vartheta\left(	\tfrac{\xi-\sigma}{2}\right)\phi\left(\tfrac{\sigma+\xi}{2},\tfrac{\xi-\sigma}{2}\right)^{-1}\right)\mathrm{d}\sigma,
\end{align*} where we used the change of variable $\sigma=2s-\xi$. In order to simplify the notations, we introduce the functions
\begin{align*}
\psi(\xi,\sigma) & :=\phi\left(\tfrac{\sigma+\xi}{2},\tfrac{\xi-\sigma}{2}\right) =\min\{\xi,\sigma\}+2R,\\
\zeta(\xi-\sigma)& :=\vartheta\left(	\tfrac{\xi-\sigma}{2}\right)=\min\left\{1,\tfrac{2R}{|\xi-\sigma|}\right\}.
\end{align*} Hence,
\begin{align*}
J_1(t,r) & =\frac{1}{2}\int_{\max\{-\xi,-R\}}^\eta \zeta(\xi-\sigma)\psi(\xi,\sigma)^{-2}\mu\left(\widetilde{C}\|u\|_{\mathrm{X}(R)}\zeta(\xi-\sigma)\psi(\xi,\sigma)^{-1}\right)\mathrm{d}\sigma \\ & \lesssim \int_{-R}^\eta \zeta(\xi-\sigma)\psi(\xi,\sigma)^{-2}\mu\left(\widetilde{C}\|u\|_{\mathrm{X}(R)}\zeta(\xi-\sigma)\psi(\xi,\sigma)^{-1}\right)\mathrm{d}\sigma =:\widetilde{J}_1(\xi,\eta).
\end{align*} Thus, recalling $\psi(\xi,\eta)=\phi(t,r)$, our goal is to prove  for any $(\xi,\eta)\in\mathbb{R}^2$ such that $\xi \geqslant -R$, $\eta\geqslant -R$ and $\xi+\eta\geqslant 0$ the estimate
\begin{align}\label{estimate J1 tilde}
\widetilde{J}_1(\xi,\eta) \lesssim  \left(\mu\left(\widetilde{C}R^{-1}\|u\|_{\mathrm{X}(R)}\right) +\int_{\widetilde{C}(t+R)^{-1}\|u\|_{\mathrm{X}(R)}}^{\widetilde{C}R^{-1}\|u\|_{\mathrm{X}(R)}} \frac{\mu(\tau)}{\tau}\, \mathrm{d}\tau\right) \psi(\xi,\eta)^{-1}.
\end{align} In order to prove \eqref{estimate J1 tilde}, we consider four different subcases, cf. \cite[Lemma 1]{Sideris=1983}.

\subsubsection*{\large Case 1: $-R\leqslant \xi\leqslant R$} In this case, since $-R\leqslant \sigma \leqslant \eta$ and $\xi\in [-R,R]$, we have that $R\leqslant \psi(\xi,\sigma), \psi(\xi,\eta)\leqslant 3R$. Therefore, to prove \eqref{estimate J1 tilde} in this case, it is sufficient to show that 
\begin{align} \label{estimate J1 tilde case1}
\widetilde{J}_1(\xi,\eta) \lesssim \mu\left(\widetilde{C}R^{-1}\|u\|_{\mathrm{X}(R)}\right) +\int_{\widetilde{C}(t+R)^{-1}\|u\|_{\mathrm{X}(R)}}^{\widetilde{C}R^{-1}\|u\|_{\mathrm{X}(R)}} \frac{\mu(\tau)}{\tau}\, \mathrm{d}\tau.
\end{align}
We remark that
\begin{align}\label{cases zeta(xi-sigma)}
\zeta(\xi-\sigma)=\begin{cases}
\dfrac{2R}{|\xi-\sigma|} & \mbox{if} \ \sigma\leqslant \xi-2R   \ \mbox{or} \  \sigma\geqslant \xi+2R, \\
1 & \mbox{if} \ \xi-2R\leqslant \sigma\leqslant\xi+2R.
\end{cases}
\end{align} Because of $\xi\in[-R,R]$, it results that $\xi-2R\leqslant -R \leqslant \sigma$ and $\xi+2R\geqslant-R$, so in case 1 from \eqref{cases zeta(xi-sigma)}  we have
\begin{align*}
\zeta(\xi-\sigma)=\begin{cases}
\dfrac{2R}{\sigma-\xi} & \mbox{if} \  \sigma\geqslant \xi+2R, \\
1 & \mbox{if} \ -R\leqslant \sigma\leqslant\xi+2R.
\end{cases}
\end{align*} 
Therefore,
\begin{align*}
\widetilde{J}_1(\xi,\eta) &\lesssim\int_{-R}^{\eta} \zeta(\xi-\sigma)\mu\left(\widetilde{C}R^{-1}\|u\|_{\mathrm{X}(R)}\zeta(\xi-\sigma)\right)\mathrm{d}\sigma \\&  \lesssim \int_{-R}^{\min\{\xi+2R,\eta\}} \mu\left(\widetilde{C}R^{-1}\|u\|_{\mathrm{X}(R)}\right)\mathrm{d}\sigma+\int^{\eta}_{\min\{\xi+2R,\eta\}} \frac{2R}{\sigma-\xi}\, \mu\left(2\widetilde{C}\|u\|_{\mathrm{X}(R)}(\sigma-\xi)^{-1}\right)\mathrm{d}\sigma \\
&  \lesssim 4R \mu\left(\widetilde{C}R^{-1}\|u\|_{\mathrm{X}(R)}\right)+\int^{\eta}_{\min\{\xi+2R,\eta\}} 2R(\sigma-\xi)^{-1}\mu\left(2\widetilde{C}\|u\|_{\mathrm{X}(R)}(\sigma-\xi)^{-1}\right)\mathrm{d}\sigma.
\end{align*} For $\eta>\xi+2R$, we consider the change of variable $\tau= 2\widetilde{C}\|u\|_{\mathrm{X}(R)}(\sigma-\xi)^{-1}$ in the last integral
\begin{align*}
\int^{\eta}_{\min\{\xi+2R,\eta\}} 2R(\sigma-\xi)^{-1}\mu\left(2\widetilde{C}\|u\|_{\mathrm{X}(R)}(\sigma-\xi)^{-1}\right)\mathrm{d}\sigma &=2R \int_{2\widetilde{C}(\eta-\xi)^{-1}\|u\|_{\mathrm{X}(R)}}^{\widetilde{C}R^{-1}\|u\|_{\mathrm{X}(R)}}\frac{\mu(\tau)}{\tau}\, \mathrm{d}\tau \\ 
& \lesssim \int_{\widetilde{C}(t+R)^{-1}\|u\|_{\mathrm{X}(R)}}^{\widetilde{C}R^{-1}\|u\|_{\mathrm{X}(R)}}\frac{\mu(\tau)}{\tau}\, \mathrm{d}\tau,
\end{align*} where in the last step we used the inequality $\eta-\xi=-2r=2|r|\leqslant 2(t+R)$. This completes the proof of \eqref{estimate J1 tilde case1}.

\subsubsection*{\large Case 2A: $\xi> R$ and $-R\leqslant \eta\leqslant \xi-2R$}

From \eqref{cases zeta(xi-sigma)} we find $\zeta(\xi-\sigma)=2R(\xi-\sigma)^{-1}$ because $\sigma\leqslant\eta\leqslant \xi-2R$. Moreover, $\psi(\xi,\eta)=\eta+2R$ and $\psi(\xi,\sigma)=\sigma+2R$ since $\sigma\leqslant\eta\leqslant \xi$. Hence,
\begin{align*}
\widetilde{J}_1(\xi,\eta) & = \int_{-R}^\eta 2R (\xi-\sigma)^{-1} (\sigma+2R)^{-2}\mu\left(2R\widetilde{C}\|u\|_{\mathrm{X}(R)}(\xi-\sigma)^{-1}(\sigma+2R)^{-1}\right)\mathrm{d}\sigma \\
& = \left(\int_{-R}^{\frac{\eta-R}{2}} +\int_{\frac{\eta-R}{2}}^\eta \right)2R (\xi-\sigma)^{-1} (\sigma+2R)^{-2}\mu\left(2R\widetilde{C}\|u\|_{\mathrm{X}(R)}(\xi-\sigma)^{-1}(\sigma+2R)^{-1}\right)\mathrm{d}\sigma \\&  =: \widetilde{J}_{1,1}(\xi,\eta)+\widetilde{J}_{1,2}(\xi,\eta).
\end{align*} For the first integral we have
\begin{align}
 \widetilde{J}_{1,1}(\xi,\eta) & =  2R\int_{-R}^{\frac{\eta-R}{2}}  (\xi-\sigma)^{-1} (\sigma+2R)^{-2}\mu\left(2R\widetilde{C}\|u\|_{\mathrm{X}(R)}(\xi-\sigma)^{-1}(\sigma+2R)^{-1}\right)\mathrm{d}\sigma \notag\\
  & \lesssim  \left(\xi-\tfrac{\eta-R}{2}\right)^{-1} \int_{-R}^{\frac{\eta-R}{2}}  (\sigma+2R)^{-2}\mu\left(\widetilde{C}\|u\|_{\mathrm{X}(R)}(\sigma+2R)^{-1}\right)\mathrm{d}\sigma \notag \\
  & \lesssim  \left(2\xi+R-\eta\right)^{-1} \mu\left(\widetilde{C}R^{-1}\|u\|_{\mathrm{X}(R)}\right) \int_{-R}^{\frac{\eta-R}{2}}  (\sigma+2R)^{-2}\mathrm{d}\sigma \notag \\
  & \lesssim  \left(\eta+2R\right)^{-1} \mu\left(\widetilde{C}R^{-1}\|u\|_{\mathrm{X}(R)}\right) = \psi(\xi,\eta)^{-1} \mu\left(\widetilde{C}R^{-1}\|u\|_{\mathrm{X}(R)}\right). \label{estimate J1,1 tilde}
\end{align}
For the second integral we obtain
\begin{align*}
\widetilde{J}_{1,2}(\xi,\eta) & =2R \int_{\frac{\eta-R}{2}}^\eta  (\xi-\sigma)^{-1} (\sigma+2R)^{-2}\mu\left(2R\widetilde{C}\|u\|_{\mathrm{X}(R)}(\xi-\sigma)^{-1}(\sigma+2R)^{-1}\right)\mathrm{d}\sigma \\
& \lesssim (\eta+3R)^{-2} \int_{\frac{\eta-R}{2}}^\eta  (\xi-\sigma)^{-1} \mu\left(2\widetilde{C}\|u\|_{\mathrm{X}(R)}(\xi-\sigma)^{-1}\right)\mathrm{d}\sigma  \\
& \lesssim (\eta+2R)^{-1} \int_{\frac{\eta-R}{2}}^\eta  (\xi-\sigma)^{-1} \mu\left(2\widetilde{C}\|u\|_{\mathrm{X}(R)}(\xi-\sigma)^{-1}\right)\mathrm{d}\sigma , 
\end{align*} where in the last inequality we used $(\eta+3R)^{-2} \leqslant (\eta+2R)^{-2} \lesssim (\eta+2R)^{-1}=\psi(\xi,\eta)^{-1} $ (recall that $\eta\geqslant -R$). Carrying out the change of variable $\tau = 2\widetilde{C}\|u\|_{\mathrm{X}(R)}(\xi-\sigma)^{-1}$, we get
\begin{align*}
\int_{\frac{\eta-R}{2}}^\eta  (\xi-\sigma)^{-1} \mu\left(2\widetilde{C}\|u\|_{\mathrm{X}(R)}(\xi-\sigma)^{-1}\right)\mathrm{d}\sigma & = \int_{4\widetilde{C}(2\xi+R-\eta)^{-1}\|u\|_{\mathrm{X}(R)}}^{2\widetilde{C}(\xi-\eta)^{-1}\|u\|_{\mathrm{X}(R)}}  \frac{\mu(\tau)}{\tau}\,\mathrm{d}\tau \\ & \leqslant \int_{\widetilde{C}(t+R)^{-1}\|u\|_{\mathrm{X}(R)}}^{\widetilde{C}R^{-1}\|u\|_{\mathrm{X}(R)}}  \frac{\mu(\tau)}{\tau}\,\mathrm{d}\tau,
\end{align*} where we employed $\eta\leqslant \xi-2R$ to increase the upper bound in the domain of integration and $2\xi+R-\eta= t+3r+R \leqslant 4(t+R)$ to decrease the lower bound in the domain of integration. So, we proved that 
\begin{align}
\widetilde{J}_{1,2}(\xi,\eta) \lesssim  \psi(\xi,\eta)^{-1} \int_{\widetilde{C}(t+R)^{-1}\|u\|_{\mathrm{X}(R)}}^{\widetilde{C}R^{-1}\|u\|_{\mathrm{X}(R)}}  \frac{\mu(\tau)}{\tau}\,\mathrm{d}\tau. \label{estimate J1,2 tilde}
\end{align} Combining \eqref{estimate J1,1 tilde} and \eqref{estimate J1,2 tilde}, we conclude the validity of \eqref{estimate J1 tilde}.

\subsubsection*{\large Case 2B: $\xi >R$ and $\xi-2R< \eta\leqslant \xi$} In this case we can split $\widetilde{J}_1(\xi,\eta)$ as follows:
\begin{align*}
\widetilde{J}_1(\xi,\eta) = \widetilde{J}_1(\xi,\xi-2R)+\int_{\xi-2R}^\eta \zeta(\xi-\sigma)\psi(\xi,\sigma)^{-2}\mu\left(\widetilde{C}\|u\|_{\mathrm{X}(R)}\zeta(\xi-\sigma)\psi(\xi,\sigma)^{-1}\right)\mathrm{d}\sigma.
\end{align*}
By the estimate in case 2A, we obtain
\begin{align}
\widetilde{J}_1(\xi,\xi-2R) & \lesssim \left(\mu\left(\widetilde{C}R^{-1}\|u\|_{\mathrm{X}(R)}\right) +\int_{\widetilde{C}(t+R)^{-1}\|u\|_{\mathrm{X}(R)}}^{\widetilde{C}R^{-1}\|u\|_{\mathrm{X}(R)}} \frac{\mu(\tau)}{\tau}\, \mathrm{d}\tau\right) \psi(\xi,\xi-2R)^{-1} \notag \\
& \lesssim \left(\mu\left(\widetilde{C}R^{-1}\|u\|_{\mathrm{X}(R)}\right) +\int_{\widetilde{C}(t+R)^{-1}\|u\|_{\mathrm{X}(R)}}^{\widetilde{C}R^{-1}\|u\|_{\mathrm{X}(R)}} \frac{\mu(\tau)}{\tau}\, \mathrm{d}\tau\right) \psi(\xi,\eta)^{-1}, \label{estimate case 2B part 1}
\end{align} 
where we used 
\begin{align*}
 \psi(\xi,\xi-2R)=\xi>\tfrac{1}{3}(\xi+2R)\geqslant \tfrac{1}{3}(\eta+2R)\approx\psi(\xi,\eta).
\end{align*} 
On the other hand, for $\xi-2R\leqslant\sigma\leqslant\eta\leqslant \xi$ we have $\psi(\xi,\sigma)=\sigma+2R$ and \eqref{cases zeta(xi-sigma)} implies that $\zeta(\xi-\sigma)=1$, consequently,
\begin{align}
& \int_{\xi-2R}^\eta \zeta(\xi-\sigma)\psi(\xi,\sigma)^{-2}\mu\left(\widetilde{C}\|u\|_{\mathrm{X}(R)}\zeta(\xi-\sigma)\psi(\xi,\sigma)^{-1}\right)\mathrm{d}\sigma  \notag \\ 
& \qquad = \int_{\xi-2R}^\eta (\sigma+2R)^{-2}\mu\left(\widetilde{C}\|u\|_{\mathrm{X}(R)}(\sigma+2R)^{-1}\right)\mathrm{d}\sigma \leqslant  \xi^{-2} \mu\left(\widetilde{C}R^{-1}\|u\|_{\mathrm{X}(R)}\right) \int_{\xi-2R}^\eta \mathrm{d}\sigma \notag \\
& \qquad =  \xi^{-2} \mu\left(\widetilde{C}R^{-1}\|u\|_{\mathrm{X}(R)}\right) (\eta-\xi+2R) \lesssim \psi(\xi,\eta)^{-1} \mu\left(\widetilde{C}R^{-1}\|u\|_{\mathrm{X}(R)}\right),  \label{estimate case 2B part 2}
\end{align} where in the last step we used $\xi\gtrsim \eta+2R= \psi(\xi,\eta)$. From \eqref{estimate case 2B part 1} and \eqref{estimate case 2B part 2} it follows \eqref{estimate J1 tilde}. 

\subsubsection*{\large Case 2C: $R<\xi <\eta$}
Similarly to case 2B, we split the integral $\widetilde{J}_1(\xi,\eta)$ into two parts
\begin{align*}
\widetilde{J}_1(\xi,\eta) = \widetilde{J}_1(\xi,\xi)+\int_{\xi}^\eta \zeta(\xi-\sigma)\psi(\xi,\sigma)^{-2}\mu\left(\widetilde{C}\|u\|_{\mathrm{X}(R)}\zeta(\xi-\sigma)\psi(\xi,\sigma)^{-1}\right)\mathrm{d}\sigma.
\end{align*} The estimate in case 2B yields
\begin{align}
\widetilde{J}_1(\xi,\xi)  \lesssim \left(\mu\left(\widetilde{C}R^{-1}\|u\|_{\mathrm{X}(R)}\right) +\int_{\widetilde{C}(t+R)^{-1}\|u\|_{\mathrm{X}(R)}}^{\widetilde{C}R^{-1}\|u\|_{\mathrm{X}(R)}} \frac{\mu(\tau)}{\tau}\, \mathrm{d}\tau\right) \psi(\xi,\eta)^{-1},  \label{estimate case 2C part 1}
\end{align} 
due to $\psi(\xi,\eta)=\xi+2R=\psi(\xi,\xi)$. From \eqref{cases zeta(xi-sigma)} we obtain for the second integral
\begin{align*}
& \int_{\xi}^\eta \zeta(\xi-\sigma)\psi(\xi,\sigma)^{-2}\mu\left(\widetilde{C}\|u\|_{\mathrm{X}(R)}\zeta(\xi-\sigma)\psi(\xi,\sigma)^{-1}\right)\mathrm{d}\sigma \\ & \quad \leqslant (\xi+2R)^{-2}\int_{\xi}^\eta \zeta(\xi-\sigma)\mu\left(\widetilde{C}R^{-1}\|u\|_{\mathrm{X}(R)}\zeta(\xi-\sigma)\right)\mathrm{d}\sigma \\ 
& \quad = \psi(\xi,\eta)^{-2}  \left[\mu\left(\widetilde{C}R^{-1}\|u\|_{\mathrm{X}(R)}\right) \int_{\xi}^{\min\{\xi+2R,\eta\}}\mathrm{d}\sigma + \int^{\eta}_{\min\{\xi+2R,\eta\}} \frac{2R}{\sigma-\xi}\, \mu\left(2\widetilde{C}\|u\|_{\mathrm{X}(R)}(\sigma-\xi)^{-1}\right)\mathrm{d}\sigma\right] \\ 
& \quad \lesssim  \psi(\xi,\eta)^{-1}  \left[\mu\left(\widetilde{C}R^{-1}\|u\|_{\mathrm{X}(R)}\right) + \int^{\eta}_{\min\{\xi+2R,\eta\}} (\sigma-\xi)^{-1} \mu\left(2\widetilde{C}\|u\|_{\mathrm{X}(R)}(\sigma-\xi)^{-1}\right)\mathrm{d}\sigma\right].
\end{align*} When $\xi+2R<\eta$ we perform the change of variable $\tau=2\widetilde{C}\|u\|_{\mathrm{X}(R)}(\sigma-\xi)^{-1}$ obtaining
\begin{align*}
\int^{\eta}_{\min\{\xi+2R,\eta\}} (\sigma-\xi)^{-1} \mu\left(2\widetilde{C}\|u\|_{\mathrm{X}(R)}(\sigma-\xi)^{-1}\right)\mathrm{d}\sigma & = \int_{2\widetilde{C}(\eta-\xi)^{-1}\|u\|_{\mathrm{X}(R)}}^{\widetilde{C}R^{-1}\|u\|_{\mathrm{X}(R)}} \frac{\mu(\tau)}{\tau}\mathrm{d}\tau \\
& \lesssim \int_{\widetilde{C}(t+R)^{-1}\|u\|_{\mathrm{X}(R)}}^{\widetilde{C}R^{-1}\|u\|_{\mathrm{X}(R)}} \frac{\mu(\tau)}{\tau}\mathrm{d}\tau,
\end{align*} where we used $\eta-\xi=-2r\leqslant 2(t+R)$.
Summarizing, we showed that
\begin{align}
& \int_{\xi}^\eta \zeta(\xi-\sigma)\psi(\xi,\sigma)^{-2}\mu\left(\widetilde{C}\|u\|_{\mathrm{X}(R)}\zeta(\xi-\sigma)\psi(\xi,\sigma)^{-1}\right)\mathrm{d}\sigma \notag \\ & \qquad\lesssim \psi(\xi,\eta)^{-1}  \left(\mu\left(\widetilde{C}R^{-1}\|u\|_{\mathrm{X}(R)}\right) + \int_{\widetilde{C}(t+R)^{-1}\|u\|_{\mathrm{X}(R)}}^{\widetilde{C}R^{-1}\|u\|_{\mathrm{X}(R)}} \frac{\mu(\tau)}{\tau}\mathrm{d}\tau\right). \label{estimate case 2C part 2}
\end{align} Hence, combining \eqref{estimate case 2C part 1} and \eqref{estimate case 2C part 2} we conclude that \eqref{estimate J1 tilde} holds even in this last case.\\

\noindent \textbf{Estimate of $J_2(t,r)$:}

From \eqref{def J2(t,r)} we have
\begin{align*}
J_2(t,r) & =\int_{\max\{0,\frac{1}{2}(t-r-R)\}}^t \vartheta(r-t+s)\phi(s,r-t+s)^{-2}\mu\left(\widetilde{C}\|u\|_{\mathrm{X}(R)}\vartheta(r-t+s)\phi(s,r-t+s)^{-1}\right)\mathrm{d}s \\
& =\frac{1}{2}\int_{\max\{-\eta,-R\}}^\xi \vartheta\left(	\tfrac{\eta-\sigma}{2}\right)\phi\left(\tfrac{\sigma+\eta}{2},\tfrac{\eta-\sigma}{2}\right)^{-2}\mu\left(\widetilde{C}\|u\|_{\mathrm{X}(R)}\vartheta\left(	\tfrac{\eta-\sigma}{2}\right)\phi\left(\tfrac{\sigma+\eta}{2},\tfrac{\eta-\sigma}{2}\right)^{-1}\right)\mathrm{d}\sigma,
\end{align*} where we carried out the change of variable $\sigma=2s-\eta$ and we used that $\vartheta, \phi(t,\cdot)$ are even functions. Clearly $J_2(t,r)\lesssim \widetilde{J}_1(\eta,\xi)$ so from \eqref{estimate J1 tilde} and $\psi(\eta,\xi)=\psi(\xi,\eta)$ we conclude that $J_2(t,r)$ can be estimated from above exactly as $J_1(t,r)$. 
\end{proof}

\subsection{Proof of Proposition \ref{Prop lower bound lifespan}} \label{Subsection lower bound lifespan 3d rad}

\hspace{5mm}We can repeat similar computations as in Subsection \ref{Subsection GESDS 3d rad}. However, since $\mu$ does not satisfy the Dini condition in this case, we can consider $\delta>0$ such that
\begin{align*}
\mu(\widetilde{C}R^{-1}\delta) &\leqslant \int_{\widetilde{C}(T+R)^{-1}\|u\|_{\mathrm{X}_T(R)}}^{\widetilde{C}R^{-1}\delta} \frac{\mu(\tau)}{\tau}\, \mathrm{d}\tau,  \\
\mu(2\widetilde{C}R^{-1}\delta) &\leqslant \int_{\widetilde{C}(T+R)^{-1}(\|u\|_{\mathrm{X}_T(R)}+\|v\|_{\mathrm{X}_T(R)})}^{2\widetilde{C}R^{-1}\delta} \frac{\mu(\tau)}{\tau}\, \mathrm{d}\tau .
\end{align*}
By using the function $\mathcal{H}$ defined in \eqref{def H function lifespan}, we can  modify suitably \eqref{norm 1/r I[u]} and \eqref{norm 1/r (I[u]-I[v])}, obtaining
\begin{align*}
\| r^{-1} I[u](t,r)\|_{\mathrm{X}_T(R)} & \leqslant 2\widetilde{C}_1 \left(\mathcal{H}(\widetilde{C}(T+R)^{-1}\|u\|_{\mathrm{X}_T(R)})-\mathcal{H}(\widetilde{C}R^{-1}\delta)\right) \|u\|_{\mathrm{X}_T(R)}^2,\\ 
 \big\| r^{-1} \big(I[u](t,r)-I[v](t,r)\big)\big\|_{\mathrm{X}_T(R)} &  \leqslant 2 \widetilde{C}_2 \left(\mathcal{H}\big(\widetilde{C}(T+R)^{-1}(\|u\|_{\mathrm{X}_T(R)}+\|v\|_{\mathrm{X}_T(R)})\big)-\mathcal{H}(2\widetilde{C}R^{-1}\delta)\right) \\ & \quad  \times\|u-v\|_{\mathrm{X}_T(R)} \left(\|u\|_{\mathrm{X}_T(R)}+\|v\|_{\mathrm{X}_T(R)}\right),
\end{align*} for any $u,v\in\mathrm{X}_T(R)$ such that $\|u\|_{\mathrm{X}_T(R)}\leqslant\delta$, $\|v\|_{\mathrm{X}_T(R)}\leqslant \delta$.

Combining the two previous inequalities with \eqref{norm u lin}, it results that $N$ is a contraction on the ball $\mathfrak{B}_\delta(\mathrm{X}_T(R))$ provided that $\delta=2\widetilde{c}_0 \|(f_0,f'_0,f''_0,f_1,f'_1)\|_{L^{\infty}([0,+\infty))} \,  \varepsilon$ and 
\begin{align*}
 \int_{\widetilde{C}(T+R)^{-1}\|u\|_{\mathrm{X}_T(R)}}^{2\widetilde{C}R^{-1}\delta} \frac{\mu(\tau)}{\tau}\, \mathrm{d}\tau   \leqslant (4\min\{\widetilde{C}_1,2\widetilde{C}_2\} \delta)^{-1}  \
\Leftrightarrow \ \mathcal{H}(\widetilde{C}(T+R)^{-1}\|u\|_{\mathrm{X}_T(R)})-\mathcal{H}(k_2 \varepsilon) \leqslant k_1 \varepsilon^{-1},
\end{align*} where 
\begin{align*}
k_1 & := (8\widetilde{c}_0 \min\{\widetilde{C}_1,2\widetilde{C}_2\} \|(f_0,f'_0,f''_0,f_1,f'_1)\|_{L^{\infty}([0,+\infty))})^{-1} , \\
 k_2 & :=4\widetilde{C} \widetilde{c}_0R^{-1} \|(f_0,f'_0,f''_0,f_1,f'_1)\|_{L^{\infty}([0,+\infty))}.
\end{align*}
Since $\mathcal{H}$ is strictly decreasing, we have that 
\begin{align*}
 \mathcal{H}(\widetilde{C}(T+R)^{-1}\|u\|_{\mathrm{X}_T(R)})\leqslant k_1 \varepsilon^{-1} +\mathcal{H}(k_2 \varepsilon) \ & \Rightarrow \ \widetilde{C}(T+R)^{-1}\|u\|_{\mathrm{X}_T(R)}\geqslant \mathcal{H}^{-1}\left(k_1 \varepsilon^{-1} +\mathcal{H}(k_2 \varepsilon) \right) \\
 & \Rightarrow \ \widetilde{C}\delta T^{-1} \geqslant \mathcal{H}^{-1}\left(k_1 \varepsilon^{-1} +\mathcal{H}(k_2 \varepsilon) \right) \\
 & \Rightarrow \ K\varepsilon T^{-1} \geqslant \mathcal{H}^{-1}\left(k_1 \varepsilon^{-1} +\mathcal{H}(k_2 \varepsilon) \right) \\
 & \Rightarrow \  T \leqslant K\varepsilon \left[\mathcal{H}^{-1}\left(k_1 \varepsilon^{-1} +\mathcal{H}(k_2 \varepsilon) \right)\right]^{-1},
\end{align*} where $K:=2\widetilde{C} \widetilde{c}_0\|(f_0,f'_0,f''_0,f_1,f'_1)\|_{L^{\infty}([0,+\infty))}$.
Thus, for any $T$ satisfying \eqref{lower bound lifespan 3d rad-T}  by Banach's fixed point argument we have a local solution in $[0,T]$. Hence, we proved the lower bound estimate for the lifespan $T_{\varepsilon}$ in \eqref{lower bound lifespan 3d rad}.

\section{Estimates for the lifespan} \label{Section upper bound lifespan}\setcounter{equation}{0}

\hspace{5mm}In this section, we unravel the upper bound estimate \eqref{Upper-Lifespan} when $\int_0^{\tau_0}\frac{\mu(\tau)}{\tau}\,\mathrm{d}\tau=+\infty$ in the cases $(1\mathrm{MoC})$ and $(2\mathrm{MoC})$ from Example \ref{Example-01}. Once more, we stress that, in the radially symmetric and 3-dimensional case, we proved in Proposition \ref{Prop lower bound lifespan} the sharpness of this estimate.
We begin, by providing a rougher estimate than the one in \eqref{Upper-Lifespan} that we will use in some cases instead of \eqref{Upper-Lifespan} to better emphasize the dominant $\varepsilon$-dependent term in the upper bound estimate for $T_\varepsilon$.

We remark that, by using de L'Hopital rule, $\ml{H}'(\tau)=-\frac{\mu(\tau)}{\tau}$ and the continuity of $\mu(\tau)$ in $\tau=0$, we have
\begin{align*}
	\lim_{\varepsilon\to 0^+} \frac{\ml{H}(c_2\varepsilon)}{c_1 \varepsilon^{-\frac{2}{n-1}}} = \lim_{\varepsilon\to 0^+} \frac{c_2\ml{H}'(c_2\varepsilon)}{-\frac{2}{n-1}c_1 \varepsilon^{-\frac{n+1}{n-1}}} = \lim_{\varepsilon\to 0^+} \frac{n-1}{2c_1}\dfrac{\dfrac{\mu(c_2\varepsilon)}{\varepsilon}}{\varepsilon^{-\frac{n+1}{n-1}}} = \lim_{\varepsilon\to 0^+} \frac{n-1}{2c_1}\mu(c_2\varepsilon)\, \varepsilon^{p_{\mathrm{Gla}}(n)-1}=0.
\end{align*}  Therefore, $\ml{H}(c_2\varepsilon)=o(\varepsilon^{-\frac{2}{n-1}})$ as $\varepsilon\to 0^+$, namely, in \eqref{Upper-Lifespan} in the argument of $\ml{H}^{-1}$ the term $\ml{H}(c_2\varepsilon)$ is dominated by the term $c_1 \varepsilon^{-\frac{2}{n-1}}$. Hence, we have  
\begin{align}\label{Upper-Lifespan 2}
	T_\varepsilon \lesssim  \varepsilon^{\frac{2}{n-1}} \left[\ml{H}^{-1}\left(2c_1 \varepsilon^{-\frac{2}{n-1}}\right)\right]^{-\frac{2}{n-1}}
\end{align} for $\varepsilon\in (0,\varepsilon_0]$ with $\varepsilon_0>0$ small enough.

Finally, we point out that, in our proof of Theorem \ref{Thm-Blow-up-sub}, the only properties of the nonlinearity $g(u_t)=|u_t|^{p_{\mathrm{Gla}}(n)}\mu(|u_t|)$ that we actually employed are the convexity and the monotonicity of $g$, the fact that $g(u_t)$ is a local nonlinearity (together with the property of finite speed of propagation for the wave equation) and the integral condition in \eqref{Assumption-Thm-Blow-up}, regardless of the fact that $\mu$ is or not a modulus of continuity. In other words, our proof is formally valid even in the cases $\mu(\tau)=\tau^{-\alpha}$ for $\alpha\in(0,\frac{2}{n-1})$ and $\mu(\tau)=1$ (clearly, these $\mu$ are not moduli of continuity), which correspond to the usual power nonlinearity of derivative-type in the sub-critical case $g(u_t)=|u_t|^p$, with $1<p<p_{\mathrm{Gla}}(n)$ and in the critical case  $g(u_t)=|u_t|^{p_{\mathrm{Gla}}(n)}$, respectively. According to these considerations, we show that \eqref{Upper-Lifespan} recovers the known sharp upper bound estimates for the lifespan of a local solution to \eqref{Wave_ut} both in the sub-critical and in the critical case.

\paragraph{\large Example (1MoC)}

In this case $\mu(\tau)=\big(\log\frac{1}{\tau}\big)^{-\gamma}$ for any $\tau\in (0,\tau_0]$, where $\gamma\in(0,1]$. Hence,
\begin{align*}
	\ml{H}(\tau)=\begin{cases}
		\displaystyle{\frac{1}{1-\gamma}\left[\left(\log\frac{1}{\tau}\right)^{1-\gamma}-\left(\log\frac{1}{\tau_0}\right)^{1-\gamma}\,\right]}&\mbox{if}\  \gamma\in(0,1),\\[1em]
		\displaystyle{\log\log\frac{1}{\tau}-\log\log\frac{1}{\tau_0}}&\mbox{if}\  \gamma=1,
	\end{cases}
\end{align*}
and
\begin{align*}
	\ml{H}^{-1}(\omega)=\begin{cases}
		\displaystyle{\exp\left\{-\left[(1-\gamma)\omega+\left(\log\frac{1}{\tau_0}\right)^{1-\gamma}\,\right]^{\frac{1}{1-\gamma}} \right\}}&\mbox{if}\  \gamma\in(0,1),\\[1.3em]
		\displaystyle{\exp\left[-\left(\log\frac{1}{\tau_0}\right)\mathrm{e}^\omega \right]}&\mbox{if}\  \gamma=1,
	\end{cases}
\end{align*}
By using \eqref{Upper-Lifespan 2} for $\gamma\in (0,1)$ and \eqref{Upper-Lifespan} for $\gamma=1$ we obtain
\begin{align*}
	T_{\varepsilon}\lesssim
	\begin{cases}
		\displaystyle{ \varepsilon^{\frac{2}{n-1}}\exp\left(c\,\varepsilon^{-\frac{2}{(n-1)(1-\gamma)}}\right)}&\mbox{if}\  \gamma\in(0,1),\\[0.8em]
		\displaystyle{ \varepsilon^{\frac{2}{n-1}}\exp\left[c_1\left(\log\frac{1}{\varepsilon}\right)\exp\left(c_2\,\varepsilon^{-\frac{2}{n-1}}\right)\right]}&\mbox{if}\  \gamma=1.
	\end{cases}
\end{align*}
We underline that if we use \eqref{Upper-Lifespan}  for $\gamma\in (0,1)$, then, as $\gamma\to 0^+$ we recover the upper bound of the sharp estimate $T_{\varepsilon} \approx \exp\left(c\,\varepsilon^{-\frac{2}{n-1}}\right)$ for the power nonlinearity with critical Glassey exponent.

\paragraph{\large Example (2MoC)}

In this case $\mu(\tau)=\prod_{j=1}^{k-1}\big(\log^{[j]}\frac{1}{\tau}\big)^{-1}\big(\log^{[k]}\frac{1}{\tau}\big)^{-\gamma}$ for any $\tau\in (0,\tau_0]$, where $k\geqslant 2$ is an integer and  $\gamma\in(0,1]$. Thus,
\begin{align*}
	\ml{H}(\tau)=\begin{cases}
		\displaystyle{\frac{1}{1-\gamma}\left[\left(\log^{[k]}\frac{1}{\tau}\right)^{1-\gamma}-\left(\log^{[k]}\frac{1}{\tau_0}\right)^{1-\gamma}\,\right]}&\mbox{if}\  \gamma\in(0,1),\\[1em]
		\displaystyle{\log^{[k+1]}\frac{1}{\tau}-\log^{[k+1]}\frac{1}{\tau_0}}&\mbox{if}\  \gamma=1,
	\end{cases}
\end{align*}
and
\begin{align*}
	\ml{H}^{-1}(\omega)=\begin{cases}
		\displaystyle{\exp\left\{-\exp^{[k-1]}\left[(1-\gamma)\omega+\left(\log^{[k]}\frac{1}{\tau_0}\right)^{1-\gamma}\,\right]^{\frac{1}{1-\gamma}} \right\}}&\mbox{if} \ \gamma\in(0,1),\\[1em]
		\displaystyle{\exp\left[-\exp^{[k]}\left(\omega+\log^{[k+1]}\frac{1}{\tau_0}\right)\right] }&\mbox{if}\  \gamma=1.
	\end{cases}
\end{align*}
Employing \eqref{Upper-Lifespan 2} for $\gamma\in (0,1)$ and \eqref{Upper-Lifespan} for $\gamma=1$ we obtain
\begin{align*}
	T_{\varepsilon}\lesssim
	\begin{cases}
		\displaystyle{ \varepsilon^{\frac{2}{n-1}}\exp\left[c_2\exp^{[k-1]}\left( c_1\,\varepsilon^{-\frac{2}{(n-1)(1-\gamma)}}\right)\right]}&\mbox{if} \ \gamma\in(0,1),\\[1em]
		\displaystyle{  \varepsilon^{\frac{2}{n-1}}\exp\left\{c_2 \exp^{[k-1]}\left[\left(\log^{[k]}\frac{1}{\varepsilon}\right)\exp\left(c_1\,\varepsilon^{-\frac{2}{n-1}}\right)\right]\right\}}&\mbox{if} \ \gamma=1.
	\end{cases}
\end{align*}

\paragraph{\large Sub-critical power nonlinearity}

Let us consider formally $\mu(\tau)=\tau^{-\alpha}$ with $\alpha\in(0,\frac{2}{n-1})$. This choice corresponds to the nonlinearity $|u_t|^{p}$ with $1<p:=p_{\mathrm{Gla}}(n)-\alpha<p_{\mathrm{Gla}}(n)$. Then,
\begin{align*}
	\ml{H}(\tau)=\frac{1}{\alpha}\left(\tau^{-\alpha}-\tau_0^{-\alpha}\right) \ \ \mbox{and}\ \ \ml{H}^{-1}(\omega)=\left( \alpha\omega+\tau_0^{-\alpha}\right)^{-\frac{1}{\alpha}} .
\end{align*}
By \eqref{Upper-Lifespan}, using that $\varepsilon^{-\frac{2}{n-1}}$ dominates $\varepsilon^{-\alpha}$,  we obtain
\begin{align*}
	T_{\varepsilon}\lesssim \varepsilon^{\frac{2}{n-1}-\frac{4}{\alpha(n-1)^2}}=\varepsilon^{-(\frac{1}{p-1}-\frac{n-1}{2})^{-1}},
\end{align*}
which coincides with the sharp upper bound estimate for the lifespan in \eqref{Sharp-Lifespan-Glassey} for the semilinear classical wave model \eqref{Wave_ut} in the sub-critical Glassey case.

\paragraph{\large Critical power nonlinearity}

Let us consider formally $\mu(\tau)=1$. In this case, the nonlinearity is exactly the power nonlinearity $|u_t|^{p_{\mathrm{Gla}}(n)}$. Thus,
\begin{align*}
	\ml{H}(\tau)=\log\tau_0-\log \tau\ \   \mbox{and}\ \    \ml{H}^{-1}(\omega)=\tau_0\,\mathrm{e}^{-\omega}.
\end{align*} By \eqref{Upper-Lifespan} we find
\begin{align*}
	T_{\varepsilon}\lesssim \exp\left(c\,\varepsilon^{-\frac{2}{n-1}}\right),
\end{align*}
which coincides with the sharp upper bound estimate for the lifespan in \eqref{Sharp-Lifespan-Glassey} for the semilinear classical wave model \eqref{Wave_ut} in the critical Glassey case.

\section*{Acknowledgments} 
 Wenhui Chen is supported in part by the National Natural Science Foundation of China (grant No. 12301270, grant No. 12171317), 2024 Basic and Applied Basic Research Topic--Young Doctor Set Sail Project (grant No. 2024A04J0016), Guangdong Basic and Applied Basic Research Foundation (grant No. 2023A1515012044). Alessandro Palmieri is a member of the Gruppo Nazionale per l’Analisi Matematica, la Probabilit\`{a} e le loro Applicazioni (GNAMPA) of the Instituto Nazionale di Alta Matematica (INdAM) and has been supported by INdAM - GNAMPA Project 2024 ``Modelli locali e non-locali con perturbazioni non-lineari'' CUP E53C23001670001 and by ERC Seeds UniBa Project ``NWEinNES'' CUP H93C23000730001. The authors thank  Michael Reissig (TU Bergakademie Freiberg) for some suggestions in the preparation of the paper.

\end{document}